\documentclass[journal]{IEEEtran}
\usepackage{tikz}
\usepackage[cmex10,intlimits]{amsmath}
\usepackage{amsfonts}
\usepackage{amssymb}
\usepackage{booktabs}
\usepackage{textcomp}

\graphicspath{figures/}
\usepackage{physics}

\providecommand{\gene}{\ensuremath{\M{g}}}

\providecommand{\setAdd}{\ensuremath{\OP{A}}}
\providecommand{\setRem}{\ensuremath{\OP{R}}}

\providecommand{\Nopt}{\ensuremath{B}}
\providecommand{\Ndof}{\ensuremath{N}}
\providecommand{\Nags}{\ensuremath{N_{\mathrm{A}}}}

\usepackage[colorlinks]{hyperref}
\hypersetup{
	colorlinks=true, 
	linkcolor=blue!70!black, 
	citecolor=red!70!black, 
	filecolor=blue!70!black, 
	urlcolor=blue!70!black, 
}

\RequirePackage{xcolor} 

\RequirePackage[cmex10,intlimits]{amsmath}
\RequirePackage{amsfonts}
\RequirePackage{amssymb}

\providecommand{\J}{\ensuremath{\mathrm{j}}}    

\providecommand{\Quot}[1]{``{#1}"}              
\providecommand{\V}[1]{\boldsymbol{#1}}         
\providecommand{\M}[1]{\mathbf{#1}}             
\providecommand{\T}[1]{\mathrm{#1}}             
\providecommand{\UV}[1]{\hat{\V{#1}}}           
\providecommand{\OP}[1]{{\mathcal{#1}}}         

\providecommand{\herm}{\mathrm{H}}

\providecommand{\ZVAC}{\ensuremath{Z_0}}           

\providecommand{\Iv}{\ensuremath{\M{I}}}
\providecommand{\Vv}{\ensuremath{\M{V}}}

\providecommand{\Xm}{\ensuremath{\M{X}}}

\providecommand{\Zm}{\ensuremath{\M{Z}}}

\providecommand{\Rmvac}{\ensuremath{\M{R}_0}}

\providecommand{\Wm}{\ensuremath{\M{W}}}

\providecommand{\Rin}{\ensuremath{R_\T{in}}}
\providecommand{\Xin}{\ensuremath{X_\T{in}}}
\providecommand{\Zin}{\ensuremath{Z_\T{in}}}

\providecommand{\Prad}{P_\T{rad}}
\providecommand{\Plost}{P_\T{lost}}

\newcommand{\ie}{\textit{i}.\textit{e}.{}} 
\newcommand{\eg}{\textit{e}.\textit{g}.{}}
\newcommand{\cf}{\textit{cf}.{}}

\RequirePackage[nolist]{acronym}
\newacro{MoM}[MoM]{method of moments}
\newacro{PEC}[PEC]{perfect electric conductor}
\newacro{RWG}[RWG]{Rao-Wilton-Glisson}
\newacro{DOF}[DOF]{\mbox{degrees-of-freedom}}
\newacro{QCQP}[QCQP]{quadratically constrained quadratic program}

\begin{document}
\title{Optimal Inverse Design Based on Memetic Algorithms -- Part 2: Examples and Properties}
\author{Miloslav~Capek, \IEEEmembership{Senior Member, IEEE}, Lukas~Jelinek, Petr~Kadlec, and Mats~Gustafsson, \IEEEmembership{Senior Member, IEEE}
\thanks{Manuscript received \today; revised \today. This work was supported by the Czech Science Foundation under project~\mbox{No.~21-19025M}.}
\thanks{M. Capek and L. Jelinek are with the Czech Technical University in Prague, Prague, Czech Republic (e-mails: \{miloslav.capek; lukas.jelinek\}@fel.cvut.cz).}
\thanks{M. Gustafsson is with Lund University, Lund, Sweden (e-mail: mats.gustafsson@eit.lth.se).}
\thanks{P. Kadlec is with the Brno University of Technology, Brno, Czech Republic (e-mail: kadlecp@feec.vutbr.cz).}
\thanks{Color versions of one or more of the figures in this paper are
available online at http://ieeexplore.ieee.org.}
}

\maketitle

\begin{abstract}
Optimal inverse design, including topology optimization and evaluation of fundamental bounds on performance, which was introduced in Part~1, 
is applied to various antenna design problems. A memetic scheme for topology optimization combines local and global techniques to accelerate convergence and maintain robustness. Method-of-moments matrices are used to evaluate objective functions and allow to determine fundamental bounds on performance. By applying the Shermann-Morrison-Woodbury identity, the repetitively performed structural update is inversion-free yet full-wave. The technique can easily be combined with additional features often required in practice, \eg{}, only a part of the structure is controllable, or evaluation of an objective function is required in a subdomain only. The memetic framework supports multi-frequency and multi-port optimization and offers many other advantages, such as an actual shape being known at every moment of the optimization. The performance of the method is assessed, including its convergence and computational cost.
\end{abstract}


\section{Introduction}

\IEEEPARstart{I}{nverse} design is a time-consuming process with no certainty regarding global minimum feasibility~\cite{Deschamps+Cabayan1972}. This is generally the case, no matter how sophisticated the method employed~\cite{BendsoeSigmund_TopologyOptimization, Ohsaki_OptimizationOfFiniteDimensionalStructures, 2018_Molesky_NatP, Angeris_etal_ANewHeuristicForPhysicalDesign, 2017_Park_Arxiv, ErentokSigmund2011, 2016_Liu_AMS, Aage_etal_TopologyOptim_Nature2017} and there is no proof of convergence towards the global minimum~\cite{BendsoeSigmund_TopologyOptimization, Ohsaki_OptimizationOfFiniteDimensionalStructures}. However, the global minimum is typically not needed in practice. A sufficiently good solution has to be found in a reasonable time. As such, a good balance between detailed local search and large-scale exploration of the solution space has to be achieved~\cite{WolpertMacready_NoFreeLunchTheoremsForOptimization}. These properties are provided by the approach introduced in Part~1~\cite{2021_capeketal_TSGAmemetics_Part1} which lays down the groundwork for this second part.

The inverse design procedure proposed in this paper relies on two equally important steps. In the first step, the fundamental bound on an optimized metric is evaluated. This provides a stopping criterion for the memetic scheme and delimits the performance of any realized device. The memetic scheme is initiated in the second step and attempts to minimize the distance from the fundamental bound. It may happen in some cases that the bound is not reachable, \eg{}, there are not enough feeders, or they are not properly placed. In such a case, the optimization can be restarted with improved initial conditions.

The memetic optimization method combines two distinct approaches -- local and global steps. The local step is based on investigating the smallest topology perturbations, \ie{}, changes in the value of an objective function if one of the \ac{DOF} is removed or added to the optimized structure. The local step was proposed in \cite{2018_ChenEtAl_AnalysisOfPartialGeometryModificationProblems} and~\cite{Capeketal_ShapeSynthesisBasedOnTopologySensitivity}, where only \ac{DOF} removals were used to detect the local minima. The approach was extended by the possibility of adding \ac{DOF} to the system in~\cite{Capeketal_InversionFreeEvaluationOfNearestNeighborsInMoM}. Satisfactory performance of the local step was confirmed on Q-factor minimization~\cite{Capeketal_InversionFreeEvaluationOfNearestNeighborsInMoM}, as well as on the minimization of reflectance of a pixel antenna~\cite{Jiang2021} or, recently, the design of surface unit cell~\cite{Wang_UnitCell}. Part~1~\cite{2021_capeketal_TSGAmemetics_Part1} merged both smallest perturbations (addition and removal of \ac{DOF}) into a unified framework and combined it with the global step.

Since a fixed discretization grid is used, the differences calculated from the change of the objective function value under the smallest topology perturbations (topology sensitivities) represent a discrete analogue to gradient over structural variables. As with gradient-based optimization schemes, these differences are used to search for a local minimum via an iterative greedy search~\cite{Cormen_etal_IntroductionToAlgorithms}.

The global step is designed to restore and maintain diversity when the local step finds a local minima. Heuristic approaches are known for their robustness~\cite{Simon_EvolutionaryOptimizationAlgorithms, OnwuboluBabu_NewOptimizationTechniquesInEngineering}, applied, \eg{}, in the form of genetic algorithms operating over locally optimal shapes. The binary nature of genetic algorithms suits the combinatorial-type optimization solved in this work. While heuristics do not, in general, perform well~\cite{Sigmund_OnTheUselessOfNongradinetApproachesInTopoOptim}, only good properties, including versatility, robustness, and easy implementation, are used here while the disadvantages (mainly computational requirements and slow convergence~\cite{Katoch2020}) are mitigated by the underlying local step. Moreover, it is shown in this paper that, in many cases, only the local step is needed to identify shapes good enough for practical purposes. The above-mentioned properties and claims are confirmed in this paper using four examples involving electrically small and medium-sized problems. Both scattering and antenna scenarios are treated. The minimizing functions are single- and multi-objective.

The paper is structured as follows. Implementation and benchmarking details are provided in Section~\ref{sec:Examples}. Section~\ref{sec:QminU} deals with the minimization of the Q-factor. Results for a rectangular plate and a spherical shell are consistently compared to fundamental bounds. The influence of discretization on precision and computational time is provided. Section~\ref{sec:TradeOff} generalizes the objective function to a multi-objective case studying the trade-off between Q-factor and input impedance (matching). A code implementing the local step for minimizing Q-factor is published as supplementary material~\cite{TSGA22_SuplMat}. An array is synthesized in Section~\ref{sec:maxSigmaS} to maximize realized gain. The different number of elements and different spacings are considered. Section~\ref{sec:maxPabs} describes how to optimize a region adjacent to a lossy chip in which the power absorbed from the incoming plane wave has to be maximized. Finally, the various aspects of the method are thoroughly discussed in Section~\ref{sec:discussion} and the paper is concluded in Section~\ref{sec:concl}.

\section{Examples -- Methodology}
\label{sec:Examples}

The properties and overall performance of the proposed optimization procedure are discussed using various examples focusing on different aspects of the method. All parts were implemented in MATLAB~\cite{matlab}, the matrix operators were evaluated in AToM~\cite{atom}, and the genetic algorithm from FOPS~\cite{fops} was utilized. The surface \ac{MoM} for good conductors was utilized with Rao-Wilton-Glisson (RWG) basis functions~\cite{RaoWiltonGlisson_ElectromagneticScatteringBySurfacesOfArbitraryShape} defined over a Delaunay triangulation~\cite{deLoeraRambauSantos_Triangulations}. The $5$-th order quadrature rule~\cite{Dunavant_HighDegreeEfficientGQR} was applied to evaluate \ac{MoM}-based integrals.

The examples with recorded computational time were evaluated on a computer with an AMD Ryzen Threadripper 1950X CPU (16 physical cores, 3.4\,GHz) with 128\,GB RAM. The extensive studies, sweeping one or more parameters, were evaluated on an RCI cluster~\cite{RCIcluster}. All parts of the code were implemented as described in~\cite{2021_capeketal_TSGAmemetics_Part1} with the exception of the GPU evaluation of an objective function, \ie{}, the local step was vectorized. The global step is calculated with the help of parallel computing.

\section{Minimum Q-factor}
\label{sec:QminU}

Minimization of Q-factor has a long history~\cite{2018_Schab_Wsto, Capek_etal_2019_OptimalPlanarElectricDipoleAntennas} and is still a topical problem. The importance of this quantity for electrically small antennas is given by its inverse proportionality to fractional bandwidth~\cite{YaghjianBest_ImpedanceBandwidthAndQOfAntennas}, a parameter suffering greatly from small electrical size~\cite{VolakisChenFujimoto_SmallAntennas}. While the determination of the lower bound to quality factor has matured~\cite{CapekGustafssonSchab_MinimizationOfAntennaQualityFactor}, the corresponding optimal shapes are known thanks only to empirical designs~\cite{Best_ElectricallySmallResonantPlanarAntennas, Fujimoto_Morishita_ModernSmallAntennas, Capek_etal_2019_OptimalPlanarElectricDipoleAntennas}. This inevitably limits their exclusive applicability to cases when they are known to be approximately optimal (\eg{}, a folded spherical helix for a spherical region~\cite{Best_LowQelectricallySmallLinearAndEllipticalPolarizedSphericalDipoleAntennas} or a meanderline for a rectangular region~\cite{Capek_etal_2019_OptimalPlanarElectricDipoleAntennas}).

Here, we demonstrate how to utilize the technique presented in Part~1 of this paper~\cite{2021_capeketal_TSGAmemetics_Part1} to get solutions close to the fundamental bounds in an automated manner. The optimization setup is chosen so that the results can be compared with empirically synthesized structures~\cite{Capek_etal_2019_OptimalPlanarElectricDipoleAntennas}, \ie{}, with the region~$\varOmega$ in the form of a rectangle as, \eg{}, described in Fig.~1 of Part~1~\cite{2021_capeketal_TSGAmemetics_Part1}.

\subsection{Inverse Design}
\label{sec:QminComp}

The rectangular region of varying mesh grid density and varying electrical size expressed in~$ka$, where~$k$ is the wavenumber and~$a$ is the radius of the smallest circumscribing sphere, is considered first. The optimization was performed with~$\Nags = 48$ agents (three times the number of physical cores of the Threadripper~1950X processor) utilized during the global step~\cite{2021_capeketal_TSGAmemetics_Part1}. The excitation is realized via a discrete delta gap feeder placed in the middle of the longer side, close to the boundary of the design region. The amplitude and phase of the excitation does not affect the result thanks to the quadratic forms used to define the Q-factor as
\begin{equation}
\begin{split}
Q (\gene) &= Q_\T{U} (\gene) + Q_\T{E} (\gene) \\
&= \dfrac{1}{2} \dfrac{\Iv^\herm \left(\gene\right) \Wm \Iv \left(\gene\right)}{\Iv^\herm \left(\gene\right) \Rmvac \Iv \left(\gene\right)} + \dfrac{1}{2} \dfrac{\left|\Iv^\herm \left(\gene\right) \Xm_0 \Iv \left(\gene\right)\right|}{\Iv^\herm \left(\gene\right) \Rmvac \Iv \left(\gene\right)}
\end{split}
\label{eq:Qfact}
\end{equation}
here decomposed into its untuned~$Q_\T{U}$ and tuning~$Q_\T{E}$ parts~\cite{2021_capeketal_TSGAmemetics_Part1}. This corresponds to the classic Q-factor definition~\cite{YaghjianBest_ImpedanceBandwidthAndQOfAntennas}
\begin{equation}
Q = \dfrac{2\omega \max \left\{ W_\T{m}, W_\T{e} \right\}}{\Prad} = \dfrac{\omega W }{\Prad} + \dfrac{\left| P_\T{react} \right| }{2\Prad},
\label{eq:QfactFromMax}
\end{equation}
where~$W, W_\T{m}, W_\T{e}$ are cycle mean values of the total, magnetic and electric stored energy, respectively, and where~$P_\T{react}, \Prad$ are the reactive and radiated powers, see Appendices~B and~C of Part~1, \cite{2021_capeketal_TSGAmemetics_Part1}. The materials allowed by the optimizer are a \ac{PEC} and vacuum.

\begin{figure}
\centering
\includegraphics[width=\columnwidth]{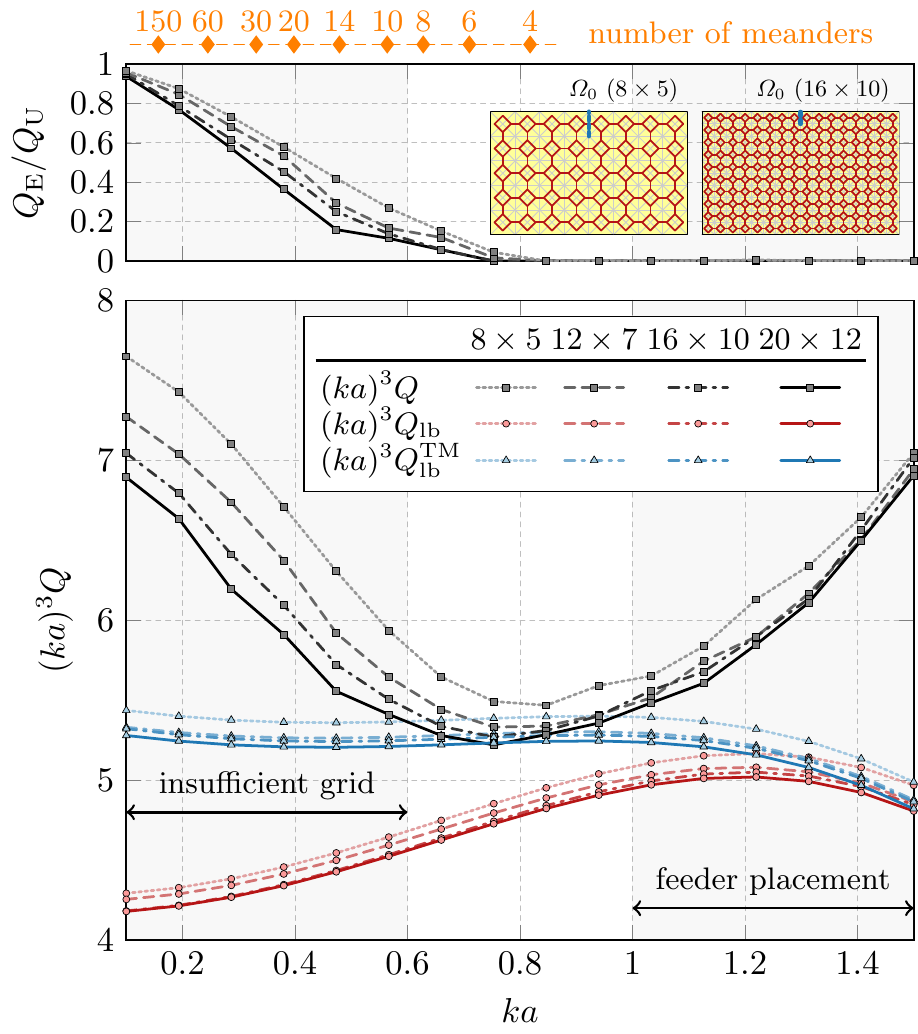}
\caption{Shape optimization (inverse design) of an antenna minimizing Q-factor. The rectangular plate is used as the region~$\varOmega$, see examples for grids $8\times 5$ and $16\times 10$ as the insets in the top-right corner. Due to cubic dependence on electrical size, Q-factor values are normalized to $(ka)^3$,~\cite{Capek_etal_2019_OptimalPlanarElectricDipoleAntennas}. The optimized performance for four mesh grids (different line styles) is compared with fundamental bounds~$Q_\T{lb}$ and~$Q_\T{lb}^\T{TM}$, \cite{CapekGustafssonSchab_MinimizationOfAntennaQualityFactor, GustafssonTayliEhrenborgEtAl_AntennaCurrentOptimizationUsingMatlabAndCVX}.  The self-resonance test $Q_\T{E}/Q_\T{U}$ is depicted in the top pane. The number of meanders which are required to construct self-resonant meanderline antenna reaching the bound~$Q_\T{lb}^\T{TM}$ and fitting the design region is indicated by diamond marks at the top. The number of meanders is inferred from the parameterization used in~\cite{Capek_etal_2019_OptimalPlanarElectricDipoleAntennas}.}
\label{figEx1bound}
\end{figure}

The algorithm was set as follows: the local step can potentially have infinitely many iterations~$I$, however, it is terminated when the relative difference~$\epsilon_\T{loc} = 10^{-7}$ between two consecutive iterations is reached. Similarly, the relative difference for the global step was set to~$\epsilon_\T{glob} = 10^{-7}$ with a maximum of $J = 250$ global iterations. 

\subsection{Comparison with the Fundamental Bounds}
\label{sec:QminPerf}

The results are compared with bounds~$Q_\T{lb}$ and~$Q_\T{lb}^\T{TM}$, shown in Fig.~\ref{figEx1bound}, and confirm the good performance of the method. The placement of the feed and reflection symmetry of the optimized region makes it possible to excite purely TM-modes~\cite{GustafssonTayliEhrenborgEtAl_AntennaCurrentOptimizationUsingMatlabAndCVX}, therefore, $Q_\T{lb}^\T{TM}$ is a tighter bound as compared to the~$Q_\T{lb}$~bound. The $Q_\T{lb}^\T{TM}$ bound is almost reached in the $ka \in [0.6,\,1]$ region. The top pane of Fig.~\ref{figEx1bound} demonstrates that the structures meeting the bound are close to self-resonance, \ie{}, $Q_\T{E}/Q_\T{U} \rightarrow 0$, following the expectation that Q-factor reaches its minimum for a self-resonant current \cite{CapekJelinek_OptimalCompositionOfModalCurrentsQ}. However, there are two regions where the bound is not met, namely, $ka < 0.6$ and $ka > 1$.

The region~$ka < 0.6$ requires finer granularity of the optimized domain than what was used in this example. This can be achieved only at the cost of increasing the number of unknowns. This is verified by the study in~\cite{Capek_etal_2019_OptimalPlanarElectricDipoleAntennas}, where the number of meanders required to construct self-resonant Q-factor-optimal meanderlines approximately increases with~$(ka)^{-2}$, see~\cite[Fig.~3]{Capek_etal_2019_OptimalPlanarElectricDipoleAntennas}. The necessary number of meanders to reach the bound is shown by orange diamond markers in Fig.~\ref{figEx1bound}. It can be verified that the bound is reachable for a mesh grid of $20\times 12$~pixels at $ka = 0.6$ where approximately $9$~meanders would be needed. This is the maximum number of meanders that this mesh grid can accommodate.

The region~$ka > 1$ is meshed densely enough, as can be seen from a comparison of the black lines in this region. However, the problem is the fixed number of feeders and their placement. It is highly probable that more than one feeder and asymmetric placements are needed to reach the lower bound on Q-factor\footnote{It should be noted that (for given material support) to every optimal current vector~$\M{I}_\T{opt}$ there exists an optimal vector of excitations~$\M{V}_\T{opt} = \M{Z} \M{I}_\T{opt}$. A study provided in~\cite{JelinekCapek_OptimalCurrentsOnArbitrarilyShapedSurfaces} further shows that even for a lower number of controllable excitations, an optimally excited current density on realized structure attempts to mimic that of the fundamental bound.}.

\begin{figure}
\centering
\includegraphics[width=\columnwidth]{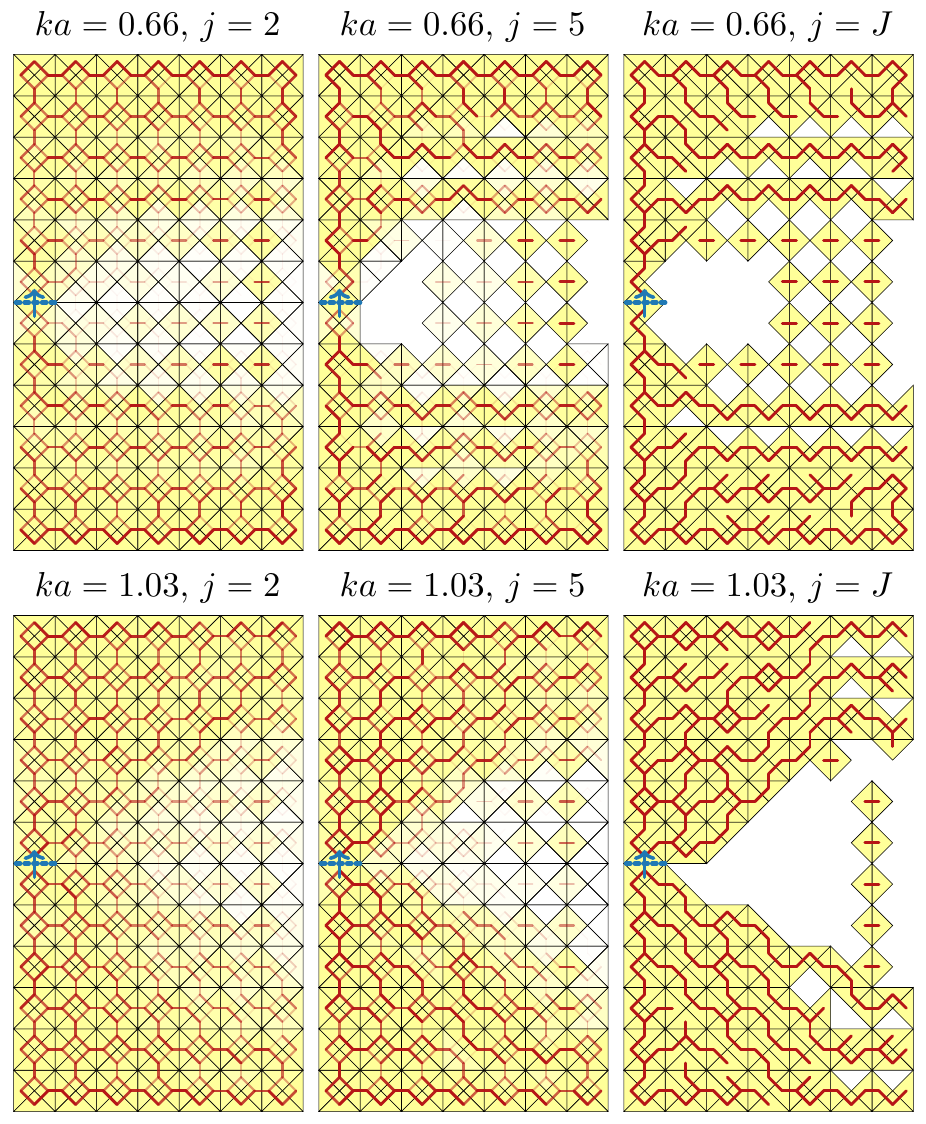}
\caption{Evolution of optimized shapes with global iteration~$j$ for a mesh grid of~$12\times 7$ square pixels. Two electrical sizes,~$ka = 0.66$ and~$ka = 1.03$ are treated \cf{}, Fig.~\ref{figEx1bound}. Only $\Nags=16$~agents were used for the global step and each picture depicts their average. A level of opacity (a value between zero and one) equals the average gene over the pool of agents of the given generation $j$, \ie{}, the more opaque the color is, the more agents have this \ac{DOF} enabled. This postprocessing technique is used only to show the diversity of the agents for a given global iteration~$j$. The last column corresponds to the final global iteration ($J = 38$ for $ka = 0.66$ and $J = 35$ for $ka = 1.03$). The solid red lines represent basis functions which are enabled, \ie{}, current paths enabled by the memetic algorithm.}
\label{figEx1agents}
\end{figure}

Two optimized shapes for~$ka = 0.66$ with~$(ka)^3 Q \approx 5.44$ and~$Q_\T{E}/Q_\T{U} \approx 0.12$, and~$ka = 1.03$ with~$(ka)^3 Q \approx 5.52$ and~$Q_\T{E}/Q_\T{U} \approx 2.3\cdot10^{-6}$ are depicted in Fig.~\ref{figEx1agents}, right column. The second and the fifth global iterations are shown as well, see the left and middle columns. In all cases, the presence of a given \ac{DOF} was evaluated as the average of all $\Nags = 48$ agents at a given $j$-th global iteration. No post-processing of the shapes was performed to show the raw results of the optimization algorithm. However, additional regularity constraints might be imposed to improve the manufacturability of the results~\cite{Capek_etal_EuCAP2021_RegularityConstraints}, \eg{}, the topology sensitivity map of the final design can be used to manually remove parts of the structure that obstruct manufacturing but have only minuscule impact on the optimized metric. It is seen that, while all agents at all iterations are always local minima (due to the performance of the local step), there is still high uncertainty about the final structure at the beginning, while the global minimum becomes uniquely determined close to the end of the optimization. Notice also, that high uncertainty at the beginning (each agent sits in different local minima) means high diversity, which is required for non-convex problems. On the other hand, the technique converges quickly, as shown in the next section.

\subsection{Computation Time}
\label{sec:QminCompTime}

The performance of the optimization technique has also to be judged in terms of computation time and convergence rate. The optimization from Fig.~\ref{figEx1bound} is repeated for~$ka = 0.5$ for varying mesh grids, see Table~\ref{tab:QuOpt}. Depending on the number of \ac{DOF}~$\Ndof$, the number of required global iterations~$J$ is shown, together with the total required computational time~$t(J)$, distance from the bound~$Q/Q_\T{lb}^\T{TM}$, and number of investigated shapes for removals ($\sum_i |\setRem(\gene_i)|$) and additions ($\sum_i |\setAdd(\gene_i)|$). The last two columns represent the number of antennas evaluated with the full-wave exact reanalysis method during the optimization. The number of optimized variables is $\Nopt = \Ndof - 1$ since one DOF is used for (fixed) delta gap feeding. It is seen that it reaches up to hundreds of millions of antenna shapes evaluated with (full-wave) \ac{MoM}.

\begin{table}[]
\caption{Total number of global iterations~$J$ used, computational time~$t$, reached Q-factor, and number of shapes evaluated for removals and additions based on the number of \ac{DOF} $\Ndof$. The values were recorded for settings from Fig.~\ref{figEx1bound} and electrical size $ka=0.5$}
\centering
\setlength{\tabcolsep}{5pt}
\begin{tabular}{ccccccc}
grid & $\Ndof$ & $J$ & $t$~(s) & $\dfrac{Q}{Q_\T{lb}^\T{TM}}$ & $\sum\limits_i |\setRem(\gene_i)|$ & $\sum\limits_i |\setAdd(\gene_i)|$ \\[2.5ex] \toprule
$8\times 5$ & $227$ & $27$ &  $24$ & $1.16$ & $8.75\cdot10^5$ & $8.78\cdot10^5$ \\
$12\times 7$ & $485$ & $49$ & $213$  & $1.13$ & $3.86\cdot10^6$ & $4.02\cdot10^6$ \\
$16\times 10$ & $934$ & $85$ &  $2190$  & $1.07$ & $1.64\cdot10^7$ & $1.84\cdot10^7$ \\
$20\times 12$ &  $1408$ & $146$ & $10400$ & $1.06$ & $3.85\cdot10^7$ & $4.36\cdot10^7$ \\
$24\times 14$ &  $1978$ & $182$ & $43236$ & $1.05$  &  $8.83\cdot10^7$ & $1.05\cdot10^8$ \\ \bottomrule
\end{tabular}
\label{tab:QuOpt}
\end{table}

The algorithmic complexity heavily depends on the implementation, computer architecture, and type of fitness function. In particular, the minimization of Q-factor requires the evaluation of quadratic forms representing stored energy and reactive power, which is computationally expensive. From this perspective, the minimization of Q-factor represents a computationally challenging example.

\begin{figure}  
\centering
\includegraphics[width=\columnwidth]{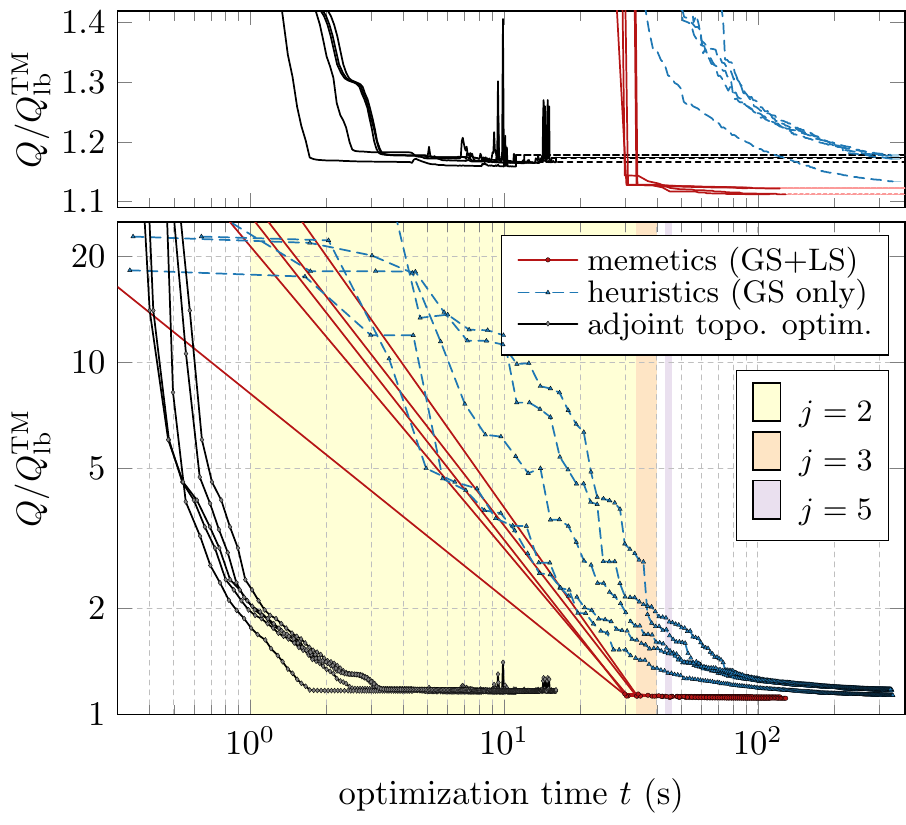}
\caption{Performance of the global step for Q-factor minimization from Fig.~\ref{figEx1bound}. The mesh grid of~$12\times 7$ pixels (dashed lines in Fig.~\ref{figEx1bound}) is used for electrical size~$ka = 0.5$. The performance is normalized to the fundamental bound~$Q_\T{lb}^\T{TM}$. Three implementations are compared: memetic combination of global and local steps proposed in this paper (red curves), the sole global step (a sole heuristic genetic algorithm, blue curves) and adjoint formulation of density-based topology optimization (black curves). Each marker represents one global iteration of the optimizer. All algorithms were run five times. The filled areas (yellow, orange, and violet) show the period when the local step of memetic scheme is performed in $j = \left\{2,3,5\right\}$ global iteration, \cf{}, Fig.~\ref{figEx1local}.}
\label{figEx1global}
\end{figure}

A graphical representation of the performance of the memetic algorithm (global and local steps together) is shown in Fig.~\ref{figEx1global} for a~$12\times 7$ mesh grid and electrical size~$ka = 0.5$, with computational time used as the $x$-axis. Only~$\Nags = 16$ agents are required for the memetic algorithm to reach very low values of the optimized metric. The remarkable convergence rate of the memetics stems from the local step. The local step, when it is applied for the first time in the second global iteration ($j = 2$), already finds a solution (after~$30$\,s) close to the final result. This occurs because of properties of the Q-factor's solution space, which makes it possible to find a reasonable solution from an arbitrary starting point, see Monte Carlo analysis in~\cite{Capeketal_InversionFreeEvaluationOfNearestNeighborsInMoM}.

\begin{figure}
\centering
\includegraphics[width=\columnwidth]{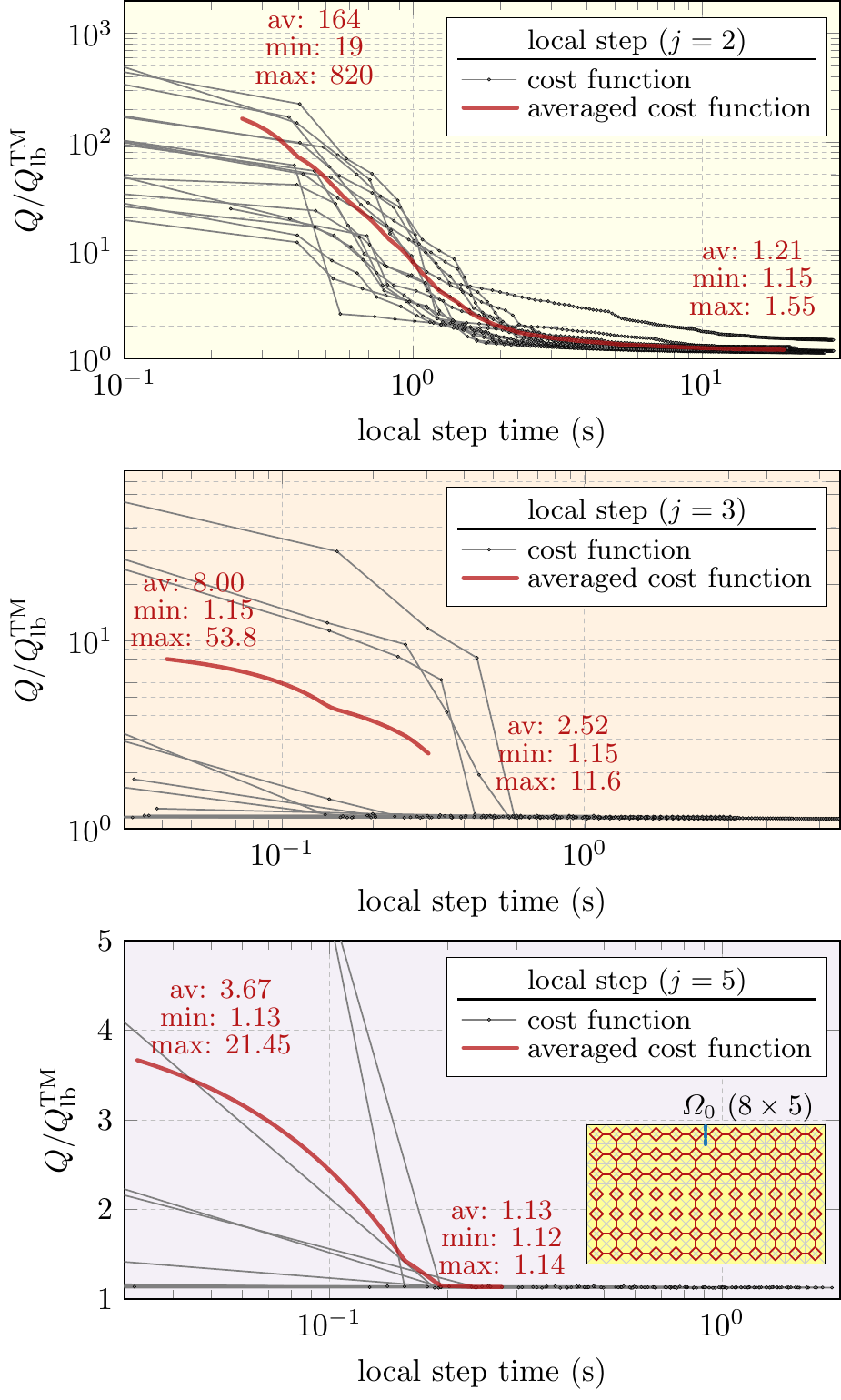}
\caption{Performance of the local step for three different global step iterations~$j = \left\{2,3,5\right\}$. A mesh grid with~$12\times7$ pixels, \ie{}, the same as in Fig.~\ref{figEx1global}, is used. The performance is normalized to the fundamental bound~$Q_\T{lb}^\T{TM}$. All~$\Nags = 16$ agents are shown in each global iteration~$j$. The computational time spent with the local updates in global iteration~$j$ is highlighted by the same background color (yellow, orange, violet) in this figure and in Fig.~\ref{figEx1global}. The thick red curve represents an average of all $\Nags$ agents (only for computational times $t$ where all curves exist). Each marker represents one local iteration of the optimizer. The minimum, maximum, and average values are shown for the averaged cost function, their beginning and end.}
\label{figEx1local}
\end{figure}

The power of the local step within the memetic scheme is emphasized in Fig.~\ref{figEx1local} where the cost function for all $\Nags = 16$~agents of memetics from Fig.~\ref{figEx1global} are depicted for global iterations $j \in \left\{2,3,5\right\}$. Iteration~$j=2$ (the top pane) is the first one when the local step is used (the first iteration only evaluates random initial seeds, see Part~I~\cite[Fig.~8]{2021_capeketal_TSGAmemetics_Part1}). Although it takes the majority of the computational time, it is also capable of decreasing the value of the objective function by three orders in magnitude, reaching, in one case,~$Q/Q_\T{lb}^\T{TM} \approx 1.15$. Then the next iteration (the middle pane) starts, on average, at a higher value of the objective function since the global operators were applied over the resulting words, consequently perturbing them and increasing the chance that a better local optimal will be found. Subsequently, the objective function is improved by one to two orders at the expense of a few seconds. Later, the fifth iteration slightly improves the local minima, almost reaching the final values of the objective function found for this grid, \cf{}, Table~\ref{tab:QuOpt}.

\subsection{Comparison with Other Optimization Schemes}

Apart from the performance of the memetic algorithm, the results obtained by sole heuristics (pure genetic algorithm, only the global step) and by density-based topology optimization~\cite{BendsoeSigmund_TopologyOptimization} are also presented in Fig.~\ref{figEx1global}. An amount of~$\Nags = 80$ agents were required for the genetic algorithm to reach values comparable to those offered by the memetic algorithm. The density-based topology optimization was using the method of moving asymptotes~\cite{Svanberg1987} with the interpolation function from~\cite{ErentokSigmund2011}. The density filter~\cite{Bruns2001} had a fixed radius of a linearly decaying cone. A continuation scheme~\cite{Wang2010} was utilized to ensure convergence to the structure consisting of either PEC or vacuum only.

The performance offered by the sheer genetic algorithm (blue lines) is inferior to the memetic scheme.

The classical topology optimization, represented by black lines in Fig.~\ref{figEx1global}, excels in speed. Consequently, it will supersede the memetic scheme when the number of DOFs is scaled up. This can be mitigated by using techniques such as surrogate modeling~\cite[Fig.~2]{Yang2019}. Density-based topology optimization is, however, not able to reach as low values of the optimized metric as compared to the memetic scheme. This is mostly given by the necessity of thresholding that converts the resulting \Quot{gray} design to a binary design made solely of PEC and vacuum. The abrupt jumps in Fig.~\ref{figEx1global} are caused by changing the properties of the filter in order to gradually remove the \Quot{gray}-scaled elements at the course of the optimization~\cite{BendsoeSigmund_TopologyOptimization}. Another drawback of classical topology optimization is the necessity to know the proper shape of the density curve and the size and shape of density filters. These are heuristic parameters not present in the memetic scheme. A notable strength of the memetic algorithm is that explicit knowledge of derivatives of the optimized metric with respect to material variables is not needed.

\subsection{Spherical Shell as a Design Region}
\label{sec:QminSphHelix}

Spherical helices fed by a delta gap source can simultaneously excite both TM and TE modes~\cite{Best_TheRadiationPropertiesOfESAsphericalHelix, Best_LowQelectricallySmallLinearAndEllipticalPolarizedSphericalDipoleAntennas}, thereby reaching the lower bound~$Q_\T{lb}$ \cite{Thal_RadiationQandGainOfTMandTEsourcesInPhaseDelayedRotatedConfigurations, CapekJelinek_OptimalCompositionOfModalCurrentsQ} and breaking the $Q_\T{lb}^\T{TM}$ bound. To construct them \textit{ad hoc} without a long study of the topic and deep empirical knowledge is, nevertheless, difficult~\cite{Thal2009}. In this section, the Q-factor without the self-resonance constraint is minimized for a spherical shell of electrical size~$ka = 0.2$, and both~$Q_\T{lb}$ and~$Q_\T{lb}^\T{TM}$ bounds\footnote{Considering electrically small region, $ka \ll 1$, the bounds for a spherical shell of radius~$a$ are $(ka)^3 Q_\T{lb} \approx 1$ and~$(ka)^3 Q_\T{lb}^\T{TM} \approx 3/2$~\cite{GustafssonCapekSchab_TradeOffBetweenAntennaEfficiencyAndQfactor}.} are utilized to judge the performance of the optimized structures.

The spherical shell was discretized into~$\Ndof = 2304$ \ac{DOF} and optimized with~$\Nags = 64$ agents and~$J = 250$ global iterations. The relative differences, both for global and local steps, were set once again to $\epsilon = 10^{-7}$. One discrete delta gap source was used (its position and orientation is irrelevant thanks to the spherical symmetry); the materials are either \ac{PEC} or vacuum.

\begin{figure}
\centering
\includegraphics[width=0.49\columnwidth]{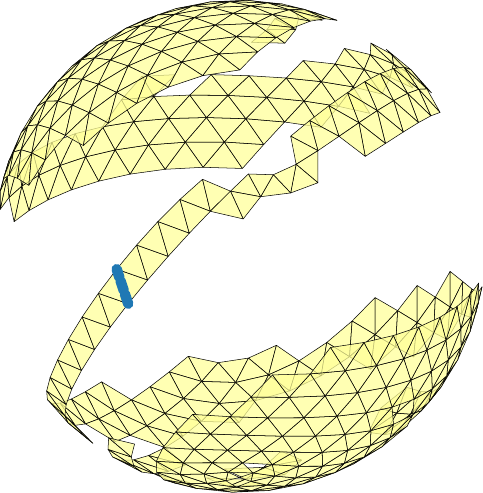}
\includegraphics[width=0.49\columnwidth]{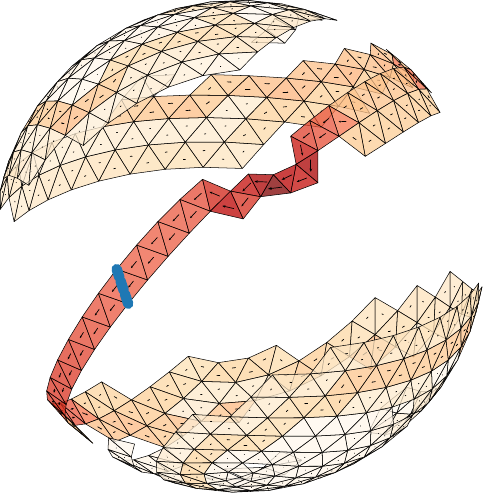}
\caption{(left) A structure optimized within a design region coinciding with a spherical shell of electrical size $ka = 0.2$, $\Ndof = 2304$. The structure is fed by a delta gap feeder depicted by the thick blue line and it resembles a spherical helix, \cf{}~\cite{Best_TheRadiationPropertiesOfESAsphericalHelix}. Normalized Q-factor reaches~$(ka)^3 Q \approx 1.05$, \ie{}, a value below the TM-bound ($Q/Q_\T{lb}^\T{TM} \approx 0.91$) and slightly above the fundamental bound ($Q/Q_\T{lb} \approx 1.27$). (right) Current density on the structure from the left pane. The absolute value of the current density is depicted by the colormap, the direction of the real part of the current is depicted by arrows. These plotting settings are used throughout the paper.}
\label{figSph1sol}
\end{figure}

\begin{figure}
\centering
\includegraphics[width=\columnwidth]{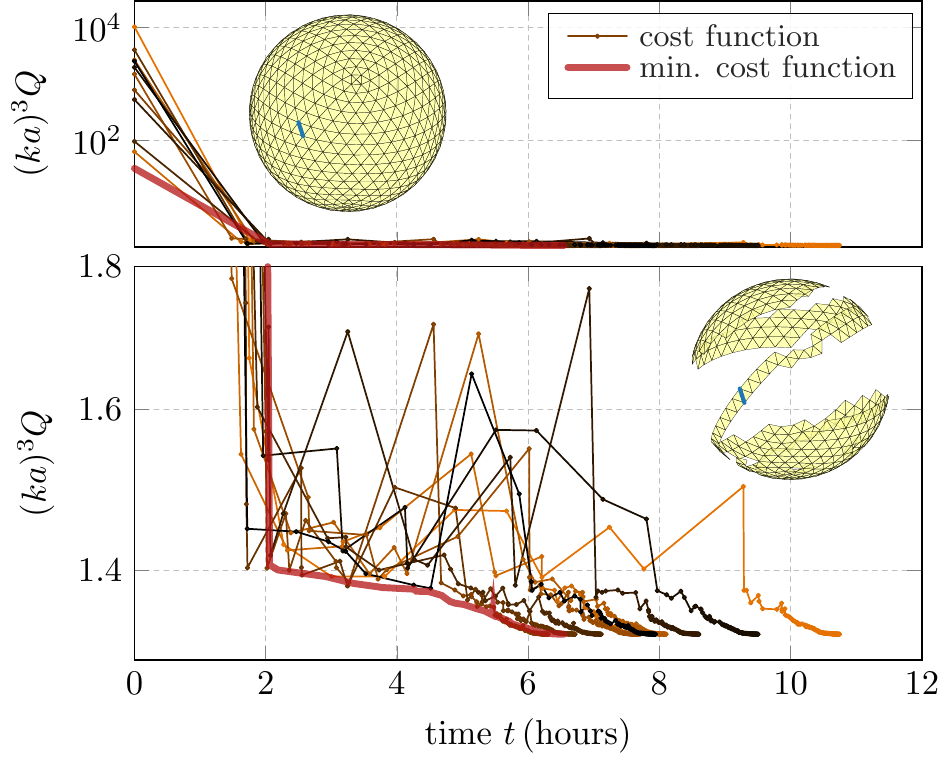}
\caption{The cost function of the first ten selected agents (out of~$64$) based on net computational time (in hours). Due to major changes in the first iterations, two panes of different $y$-axis limits are shown. The cost functions do not decrease monotonically during the global step since heuristics perturb the reached local minima. The thick red curve shows the minimum reached in each iteration of global step~$j$. The inset in the top-left corner shows the design region fully populated with \ac{DOF} at the beginning of the optimization, and the bottom right inset shows the optimized structure, \cf{}, \ref{figSph1sol}.}
\label{figSph1time}
\end{figure}

The optimization ran for $46$~hours and the resulting structure is shown in Fig.~\ref{figSph1sol}. While the computational time might seem enormous, in total $6.24\cdot 10^8$~antennas were evaluated via the exact reanalysis local step (within a solution space containing $2^{2303}\approx 1.87\cdot 10^{693}$ antenna variants). The crucial role of the local step is, once again, visible in Fig.~\ref{figSph1time} where the first iteration with a greedy search decreases the objective function from $(ka)^3 Q \approx 10^3$ to $(ka)^3 Q \approx 1.4$. Good performance in the Q-factor, similar to what was reported in~\cite{Best_LowQelectricallySmallLinearAndEllipticalPolarizedSphericalDipoleAntennas}, where the antenna was designed empirically, is observed, \ie{}, $Q/Q_\T{lb} \approx 1.27$ and $Q/Q_\T{lb} \approx 0.91$. The helix has a slope optimized for an optimal combination of TM and TE modes~\cite{CapekJelinek_OptimalCompositionOfModalCurrentsQ} and several turns to reach self-resonance. The input impedance is, however, only~$Z_\T{in} \approx 2.14\,\Omega$, which seriously limits the practical usage of such an antenna. It can be increased by prioritizing low reflectance instead of Q-factor in a multi-objective optimization as shown in the next section.

\section{Trade-off Between Q-factor and Input Impedance}
\label{sec:TradeOff}

It was demonstrated in the previous example that more than one metric is of concern in antenna design. For example, to minimize the Q-factor and to reasonably match the antenna to a given input impedance~$\Zin = \Rin + \J\Xin$ simultaneously requires composite objective functions. It can, for example, be of the following form
\begin{equation}
f\left(\gene\right) = \dfrac{Q\left(\gene\right)}{Q_\T{lb}} \left( 1 + \zeta \left\vert \varGamma \left( \gene, \Zin^0\right) \right\vert^2\right),
\label{eq:QminZin}
\end{equation}
where~$\varGamma \left( \gene, \Zin^0\right)$ is the reflection coefficient
\begin{equation}
\varGamma \left( \gene, \Zin^0\right) = \dfrac{\Zin (\gene) - \Zin^0}{\Zin (\gene) + \Zin^0}
\label{eq:Gamma}
\end{equation}
of a structure represented by a word~$\gene$, evaluated for characteristic impedance~$\Zin^0$~\cite{Pozar_MicrowaveEngineering}, and $\zeta$ is a weighting coefficient. For~$\zeta = 0$ the optimization is the same as in Section~\ref{sec:QminU}. Any value~$\zeta > 0$ takes into account the matching as well. This results in a Pareto-type optimization with a scalarized objective function~\cite{Ehrgott_MulticriteriaOptimization}, where the parameter~$\zeta$ sweeps over the feasible set.

\begin{figure}
\centering
\includegraphics[width=\columnwidth]{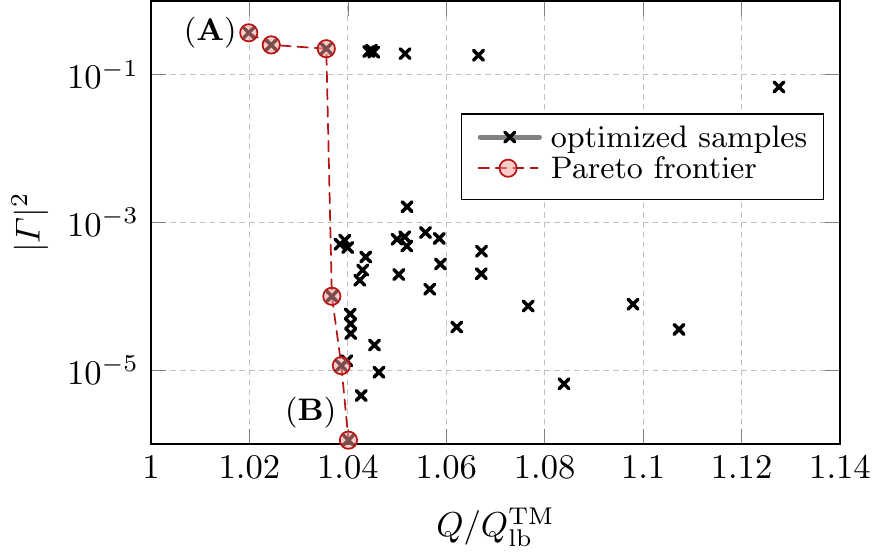}
\caption{Trade-off between low Q-factor~$Q/Q_\T{lb}^{\T{TM}}$ and matching to $Z_\T{in}^0 = 50\,\Omega$. A \ac{PEC} rectangular plate of $12\times 7$~pixels, discretized into $\Ndof = 485$~\ac{DOF} was used. The electrical size is $ka = 0.7$. The delta gap feeder is placed in the top middle. In total, $41$~samples were evaluated with weight $\zeta$ in~\eqref{eq:QminZin} set from 0 to 5 with an equidistant step. The Pareto frontier is highlighted by the red dashed line. Pareto-optimal solutions are shown as red marks. The two extreme cases, the one with the lowest Q-factor and the one well-matched to $Z_\T{in}^0$ are shown in Fig.~\ref{figQminZinTradeOffStructure}.}
\label{figQminZinTradeOff}
\end{figure}

\begin{figure}
\centering
\includegraphics[width=\columnwidth]{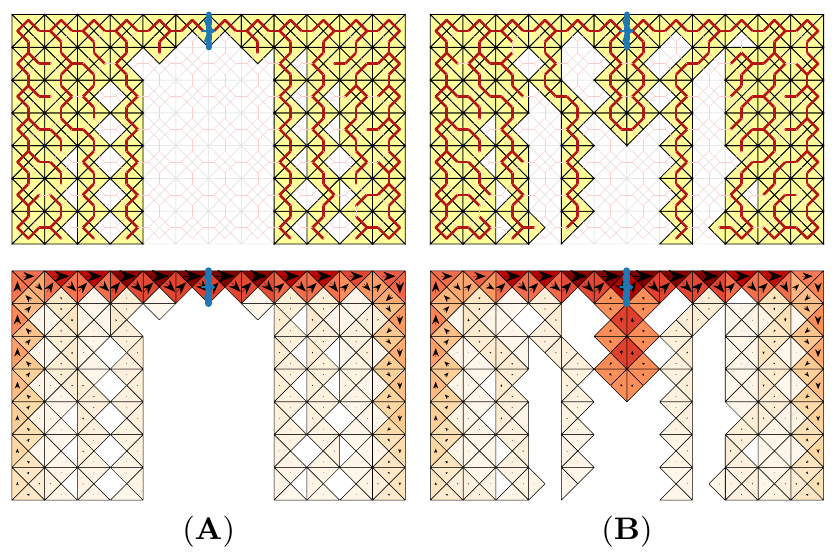}
\caption{The optimal shapes (top) and surface current densities (bottom) for two members of the Pareto frontier in Fig.~\ref{figQminZinTradeOff}. Enabled \ac{DOF} and material used are highlighted. The structures are fed by a delta gap feeder, depicted by the thick blue line. The red lines denote enabled \ac{DOF}. The (A) case minimizes Q-factor and resembles a meanderline. The (B) case is matched, therefore, there is a stub close to the feeder.}
\label{figQminZinTradeOffStructure}
\end{figure}

The PEC rectangular plate of electrical size~$ka = 0.7$ from Section~\ref{sec:QminPerf} is considered first. Self-resonance is reachable for this electrical size, and only the real part of the input impedance has to be optimized. The desired characteristic impedance is set to~$\Zin^0 = 50\,\Omega$ and weight~$\zeta$ is swept from~$\zeta = 0$ to~$\zeta = 5$ in $41$~equidistant samples. The results are depicted in Fig.~\ref{figQminZinTradeOff}. The Pareto-optimal solutions~\cite{1978CohonMultiobjectiveProgrammingAndPlanning} are highlighted by red circles and interconnected to form an approximate Pareto frontier. It is seen that the Q-factor, in terms of $Q/Q_\T{lb}^\T{TM}$, can be minimized to $Q/Q_\T{lb}^\T{TM} \approx 1.02$ but at that point the antenna is not matched. Conversely, a slight increase in~$Q/Q_\T{lb}^\T{TM}$ leads to excellent matching. The two most distinct solutions, marked by (A) and (B) in~Fig.~\ref{figQminZinTradeOff} are shown in Fig.~\ref{figQminZinTradeOffStructure} in terms of optimal structure (top) and current density (bottom). The left structure in Fig.~\ref{figQminZinTradeOff}, denoted as (A), performs best in terms of Q-factor. The right structure performs well in Q-factor and is matched to $50\,\Omega$. The visual comparison reveals that the parallel stub was created in structure (B) to match the input impedance. This technique~\cite{Fujimoto_Morishita_ModernSmallAntennas} is also used in practice~\cite[Fig.~9]{Best_LowQelectricallySmallLinearAndEllipticalPolarizedSphericalDipoleAntennas}.

\begin{figure}
\centering
\includegraphics[width=0.49\columnwidth]{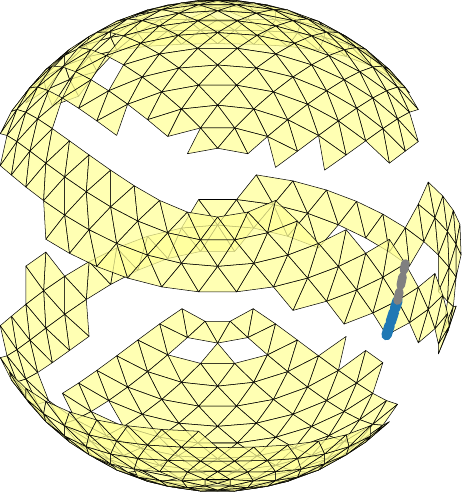}
\includegraphics[width=0.49\columnwidth]{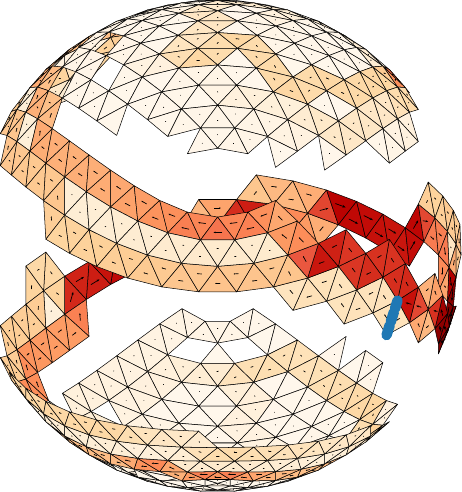}
\caption{(left) A structure optimized within a design region coinciding with a spherical shell of electrical size $ka = 0.2$, $\Ndof = 2304$. The structure is fed by a delta gap depicted by the thick blue line. The optimization considers Q-factor minimization and matching to $\Zin^0 = 50\,\Omega$. (right) Current density on the structure from the left pane. Similar to Fig.~\ref{figSph1sol}, normalized Q-factor reaches~$(ka)^3Q \approx 1.05$, \ie{}, $Q/Q_\T{lb} \approx 1.27$, and $Q/Q_\T{lb}^\T{TM} \approx 0.91$, but the reached input impedance is $\Zin \approx (51.5 - 2.39\J)\,\Omega$ ($\vert\varGamma\vert^2 \approx 7.83\cdot10^{-4}$).}
\label{figQminZinSphereStruc}
\end{figure}

The same formulation~\eqref{eq:QminZin} is applied to an optimization of the spherical shell presented in Section~\ref{sec:QminSphHelix} with an attempt to modify the structures from Fig.~\ref{figSph1sol} to be matched to~$\Zin^0 = 50\,\Omega$. For this purpose, weight $\zeta$ in \eqref{eq:QminZin} is set to $\zeta = 1$. The optimization found a minimum with $Q/Q_\T{lb} \approx 1.27$ and $Q/Q_\T{lb}^\T{TM} \approx 0.91$, \ie{}, the same as before, and $\Zin \approx (51.5 - \J 2.39)\,\Omega$, \ie{}, $|\varGamma|^2 \approx 7.72 \cdot 10^{-4}$, which is sufficient for a majority of applications~\cite{Balanis_ModernAntennaHandbook}. The optimal structure is shown in Fig.~\ref{figQminZinSphereStruc}. When the optimal shapes and currents for an unmatched antenna in Fig.~\ref{figSph1sol} and a matched antenna in Fig.~\ref{figQminZinSphereStruc} are compared, a similar difference, as in Fig.~\ref{figQminZinTradeOffStructure}, is observed, \ie{}, the optimal structure is only slightly perturbed in the vicinity of the feeder, using a stub-matching technique to divide the current flowing through the feeding port. Another improvement is a wider metallic strip used to curve the helix antenna (at least two \ac{DOF} are enabled per width of the strip).

Considering Q-factor, this section shows that there are many local minima close to the global minimum. However, other antenna metrics, such as input impedance, have different values in these local minima. In other words, Q-factor minimization can always be extended to a multi-objective case with a constraint on another metric without suffering a major increase in Q-factor. This is consistent with empirical evidence known in the literature~\cite{Best_LowQelectricallySmallLinearAndEllipticalPolarizedSphericalDipoleAntennas}.

\section{Maximal Realized Gain of an Array}
\label{sec:maxSigmaS}

Electrically larger structures are considered in this example. For this reason, the structure is preselected (an antenna array) and optimized further. A thin-strip linear antenna array consisting of $N_\T{dip}$ dipoles, operating at $1\,$GHz, and made of copper, $\sigma = 5.96\cdot10^7\,\T{Sm}^{-1}$, is generated. The thin-sheet model is utilized\footnote{It is assumed that the thickness of the structure is much larger than the penetration depth. The current is modeled as exponentially decaying from the surface of infinite half-space~\cite{Jackson_ClassicalElectrodynamics}.} to determine the surface resistivity as
\begin{equation}
R_\T{s} = \sqrt{\dfrac{k \ZVAC}{2 \sigma}},
\label{eq:surfRes}
\end{equation}
where~$\ZVAC$ is the impedance of free space. The initial length of the dipoles is~$\ell/\lambda = 0.55$ and is intentionally longer than~$\ell = \lambda/2$ to give the algorithm a chance to adopt it freely. The width of the strips is~$\ell/60$, and the separation distance is~$d$. The first element from the left always acts as a passive reflector, while the second dipole from the left is always fed by the delta gap feeder in the middle. Other dipoles act as passive directors.

Maximal realized gain~$G_\T{r}$ is chosen as a figure of merit. Since the optimization algorithm performs the minimization, the sign is switched to minus as
\begin{equation}
f\left(\gene\right) = -G_\T{r} \left(\gene\right) = - G (\gene, \UV{d}, \UV{e}) \left( 1 - \left\vert \varGamma \left( \gene, \Zin^0\right) \right\vert^2\right)
\label{eq:GmaxZin}
\end{equation}
where $G (\gene, \UV{d}, \UV{e})$ is antenna gain
\begin{equation}
G (\gene, \UV{d}, \UV{e}) = \dfrac{U (\gene, \UV{d}, \UV{e})}{\Prad (\gene) + \Plost (\gene)} = \dfrac{U (\gene, \UV{d}, \UV{e})}{P_\T{tot}(\gene)}
\label{eq:gain}
\end{equation}
and $\varGamma \left( \gene, \Zin^0\right)$ is reflection coefficient~\eqref{eq:Gamma}. From an implementation point of view, both radiation intensity and total power are one order cheaper to evaluate than Q-factor. As the strip (1D) structures are optimized, we can, therefore, expect relatively fast convergence.

\begin{figure}
\centering
\includegraphics[width=\columnwidth]{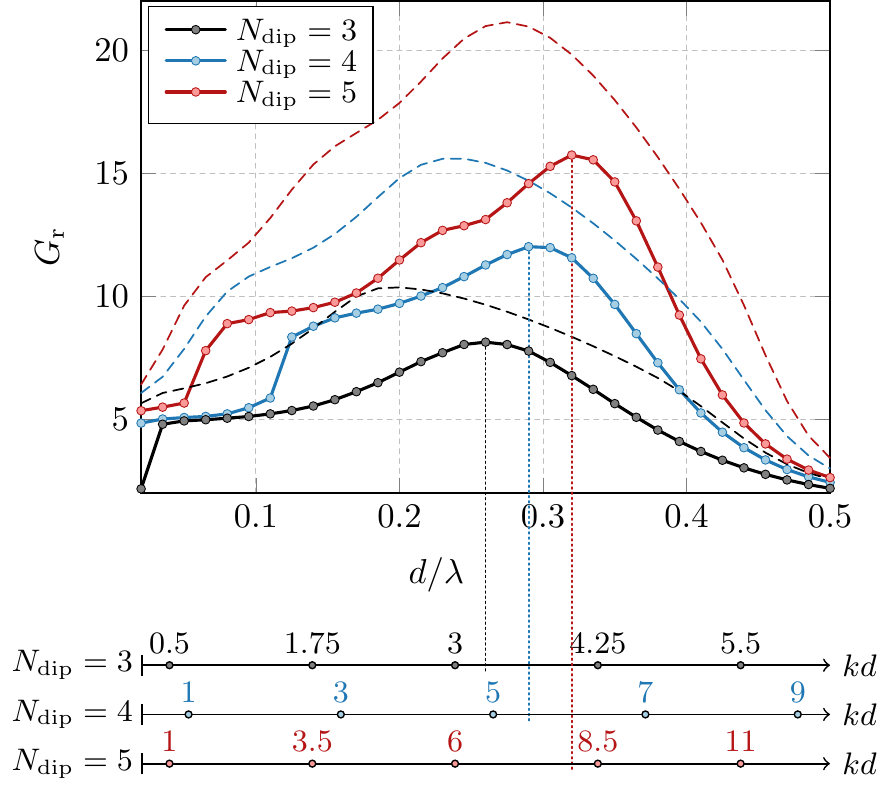}
\caption{Shape optimization of an antenna array consisting of $N_\T{dip} = \left\{3,4,5\right\}$ thin-strip dipoles, made of copper ($\sigma = 5.96\cdot10^7\,\T{Sm}^{-1}$) and operating at $1$\,GHz. The initial length of the optimized dipoles was~$\ell/\lambda = 0.55$ and the width~$w = \ell/60$. The optimization was repeated for various separation distances denoted as $d/\lambda$ to maximize realized gain. The reference impedance was set to $\Zin^0 = 50\,\Omega$, polarization pointed along the dipoles, and gain was measured in the end-fire direction. Excitation was performed by a delta gap feeder placed in the middle of the second dipole from the left, see Fig.~\ref{figGrealizedStructure}. The vertical dotted lines point to the maximal realized gain found for a various number of dipoles and denote the corresponding electrical length~$kd$. For convenience, the results are compared with fundamental bounds with prescribed input impedance (dashed lines). For one optimized sample, see Fig.~\ref{figGrealizedStructure}.}
\label{figGrealizedArray}
\end{figure}

The minimization of~\eqref{eq:GmaxZin} was repeated for $N_\T{dip} = \left\{3,4,5\right\}$ dipoles with separation distance selected from~$d/\lambda=0.02$ to~$d/\lambda = 0.5$, see the solid lines in Fig.~\ref{figGrealizedArray}. To be as thorough as possible, the fundamental bound on the realized gain with $Z_\T{in}^0 = 50\,\Omega$ was also evaluated (see the dashed lines). The procedure is detailed in Appendix~\ref{app:AppendixA}. 

There are many conclusions that can be drawn from Fig.~\ref{figGrealizedArray}. First, there is a \Quot{sweet spot} in separation distance~$d/\lambda$ (or in~$kd$ which denotes the electrical length of an array) and this optimal separation distance is different for optimized structures and for fundamental bounds (fundamental bounds reach maxima for lower distance~$kd$). Second, a higher number of array elements leads to significantly higher realized gain, see Table~\ref{tab:realizedGain} for numerical comparison. This is consistent with array theory~\cite{Balanis_Wiley_2005}. Finally, the maxima given by the fundamental bounds are closely followed by realized designs. As the desired input impedance was prescribed for the bounds, see Appendix~\ref{app:AppendixA}, this might indicate that all the realized arrays were sufficiently well matched and that the antenna gain for both bounds and designs are comparable. This has been verified by an inspection of the optimization data.

\begin{table}[]
\caption{Comparison of realized gain maxima found via memetic shape optimization and the fundamental bounds depending on the number of dipoles used. Both separation distance~$d/\lambda$ and electrical length of the array~$k d$ are shown}
\centering
\begin{tabular}{ccccccc}
& \multicolumn{3}{c}{memetics} & \multicolumn{3}{c}{fundamental bound} \\[1ex]
$N_\T{dip}$ & $d/\lambda$ & $kd$ & $G_\T{r}$ & $d/\lambda$ & $kd$ & $G_\T{r}$ \\ \toprule
$3$ & $0.260$ & $3.27$ & $\mathbf{8.13}$ & $0.200$ & $2.51$ & $\mathbf{10.4}$ \\
$4$ & $0.290$ & $5.47$ & $\mathbf{12.0}$ & $0.245$ & $4.62$ & $\mathbf{15.6}$ \\
$5$ & $0.320$ & $8.04$ & $\mathbf{15.7}$ & $0.275$ & $6.91$ & $\mathbf{21.1}$ \\ \bottomrule
\end{tabular}
\label{tab:realizedGain}
\end{table}

\begin{figure}
\centering
\includegraphics[width=\columnwidth]{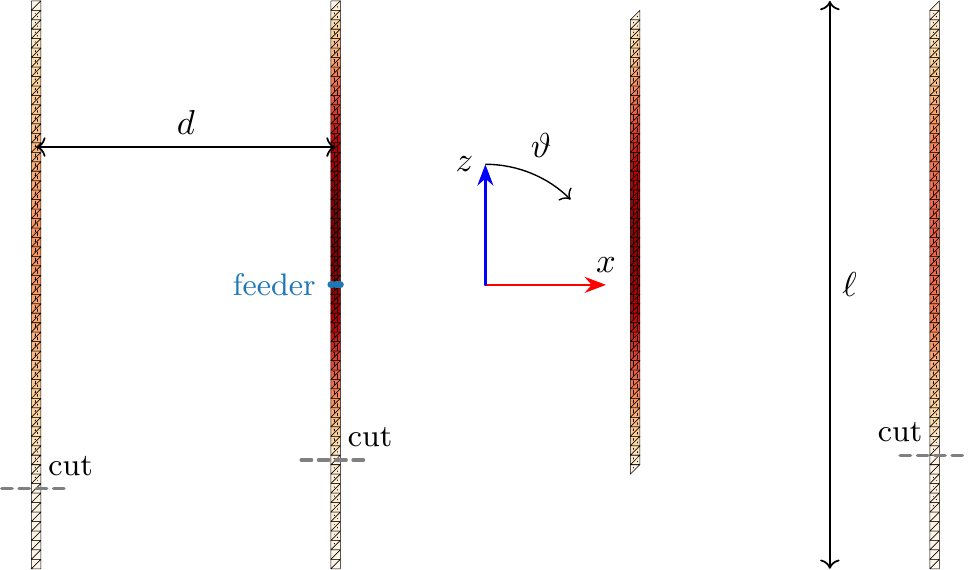}
\caption{Optimized array of four thin-strip dipole elements made of copper for maximal realized gain at $1$\,GHz. The separation distance was fixed to $d/\lambda = 0.29$, \ie{}, for a distance where the highest realized gain is reached, \cf{}, Fig.~\ref{figGrealizedArray}. The array is fed in the middle of the second dipole to the left. The reference impedance is~$\Zin^0=50\,\Omega$, polarization~$\UV{e} = \UV{z}$, and direction of main lobe~$\UV{d} = \UV{x}$. The initial electric length~$\ell/\lambda = 0.55$ is modified by cutting three dipoles. One dipole is modified by the complete removal of the material.}
\label{figGrealizedStructure}
\end{figure}

As an example, the optimal design for~$N_\T{dip} = 4$ dipoles and a separation distance chosen for the highest realized gain ($d/\lambda = 0.29$) is shown in Fig.~\ref{figGrealizedStructure}. It is seen that the structure was significantly modified. A driven element was split by the removal of one \ac{DOF}. This is indicated by the gray dashed line and a tag \Quot{cut} in Fig.~\ref{figGrealizedStructure}. The same result occurred for the first and fourth elements (counted from the left). The third element was modified by the complete removal of material from its lower part. It is obvious that these modifications cannot be achieved by a simple parametric sweep to optimize the overall length of the dipoles. Such an approach would remove all material below the \Quot{cut} label, reducing the performance from $G_\T{r} \approx 12.0$ to $G_\T{r} \approx 10.9$. The resulting structure is not symmetrical, even though the initial problem was. This is partly because the underlying discretization grid is not symmetrical and partly because the optimization problem is non-convex.

\section{Maximum Absorption in a Given Region}
\label{sec:maxPabs}

\begin{figure}
\centering
\includegraphics[width=\columnwidth]{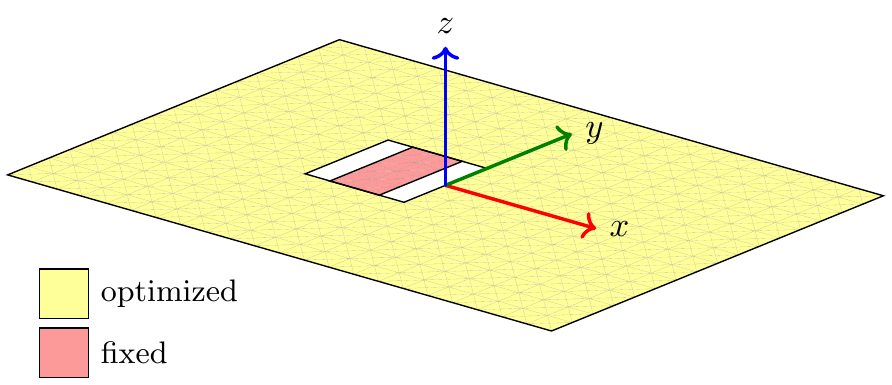}
\caption{A schematic for maximization of power absorbed in the prescribed region made of carbon ($\sigma = 1\cdot10^4\,\T{Sm}^{-1}$). The operation frequency is $f =3\,$GHz, the shorter side has length~$\lambda/4$. The amount of power absorbed from the plane wave impinging from the~$\UV{d} = -\UV{x}$ direction with~$\UV{e} = \UV{y}$ polarization is evaluated only in this region denoted by the pink color. The optimization domain is made of copper ($\sigma = 5.96\cdot10^7\,\T{Sm}^{-1}$) and denoted by the yellow color.}
\label{figmaxPabsSchematics}
\end{figure}

The last example deals with the maximization of power absorbed in a given region, \eg{}, an optimal receiving structure for an RFID chip. The optimization domain is, in this case, only a part of the entire structure, see Fig.~\ref{figmaxPabsSchematics} for a schematic layout. The controllable part, which is a subject of the optimization (highlighted by the yellow color) is made of copper, $\sigma = 5.96\cdot10^7\,\T{Sm}^{-1}$. The chip (pink color) is made of carbon, $\sigma = 1\cdot10^4\,\T{Sm}^{-1}$. The thin-sheet resistivity model~\eqref{eq:surfRes} is applied. The operational frequency is~$3\,$GHz, the shorter side has length~$\lambda/4$ and the longer~$11 \lambda/32$, and the incident plane wave is impinging from the end-fire direction, \ie{}, along the~$-\UV{x}$ direction with~$\UV{e} = \UV{y}$ polarization, see Fig.~\ref{figmaxPabsSchematics}.

The absorbed power~$P_\T{lost}$ is evaluated from
\begin{equation}
P_\T{lost} = \dfrac{1}{2} \Iv^\herm \left( \M{D}_\T{chip}^\herm \M{R}_\rho \M{D}_\T{chip} \right) \Iv,
\label{eq:Pabs}
\end{equation}
where $\Iv$ is the current on the optimized structure, $\M{R}_\rho$ is the lossy matrix, and matrix~$\M{D}_\T{chip}$ is an indexation matrix having zeros everywhere except for diagonal positions corresponding to the \ac{DOF} lying in the fixed region, see Fig.~\ref{figmaxPabsSchematics}. To ease the computational burden, the matrix~$\M{D}_\T{chip}$ is used to index out only the relevant entries required to evaluate~\eqref{eq:Pabs}.

\begin{figure}
\centering
\includegraphics[]{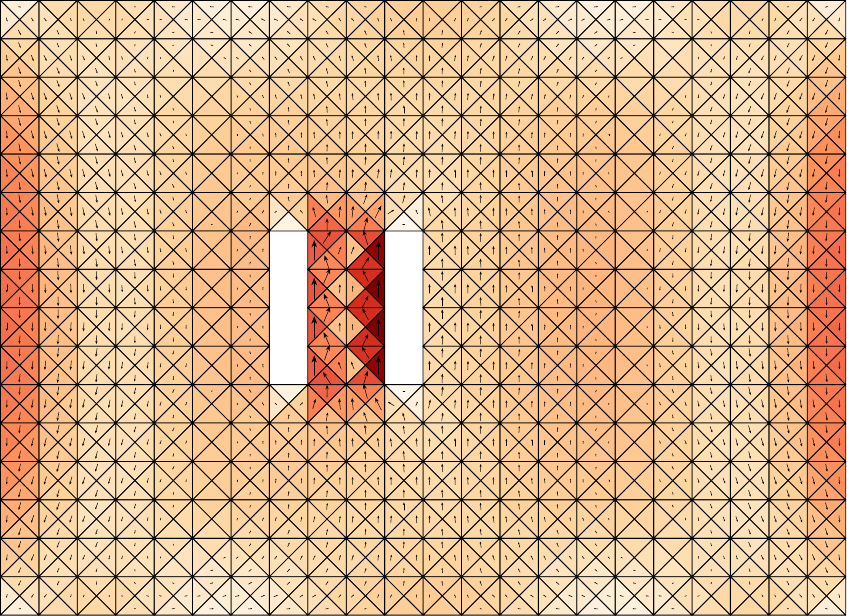}
\caption{Current density representing the optimal performance for maximal power absorbed~$P_\T{lost}$ in the uncontrollable region, depicted by the pink color in Fig.~\ref{figmaxPabsSchematics}.}
\label{figmaxPabsBound}
\end{figure}

The upper bound for this optimization problem is $P_\T{lost}^\T{ub} \approx 11.2$\,\textmu W and was found using a procedure specified in Appendix~\ref{app:AppendixB}. The optimal current impressed in vacuum is depicted in Fig.~\ref{figmaxPabsBound}. Even though this current does not fulfill $\Zm \Iv = \Vv$ with a plane wave excitation represented by a vector of expansion coefficients~$\Vv$, it is seen that it is a smooth function with maxima along the shorter sides and in the area of the chip. It must be noted that there is no guarantee that the bound is tight in this case.

\begin{figure}
\centering
\includegraphics[width=\columnwidth]{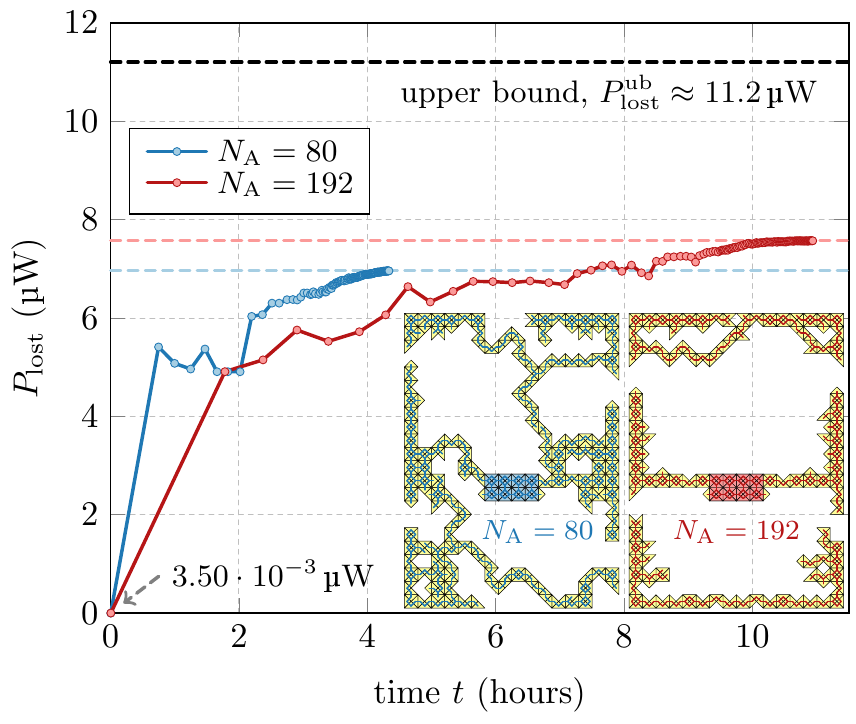}
\caption{Cost function for maximization of absorbed power~$P_\T{lost}$. Optimization was performed for $\Nags = 80$ (blue line) and $\Nags = 192$ (red line) agents. The optimal structures for both cases are shown in the bottom right corner. The structure for $\Nags = 192$ is shown in detail in Fig.~\ref{figmaxPabsCurrent}.}
\label{figmaxPabsComparison}
\end{figure}

The optimization was performed twice using the same settings and only varying the number of agents~$\Nags$ used for the global step. As in all previous examples, the first application of the local step led to an immense improvement of the objective function value, here, of absorbed power~$P_\T{lost}$ from approximately $3.50\cdot10^{-3}$\,\textmu W to $5.10$\,\textmu W, \ie{}, by more than three orders in magnitude. For $\Nags = 80$ agents (five times the number of cores of the Threadripper~1950X), the optimization ran for approximately $4.3$~hours with maximum $P_\T{lost} \approx 6.96$\,\textmu W, see Fig.~\ref{figmaxPabsComparison}. Increasing the number of agents to~$\Nags = 192$ led to maximum~$P_\T{lost} \approx 7.75$\,\textmu W found in $10.9$~hours. This value reaches $67.6\%$ of the upper bound realized by the current shown in Fig.~\ref{figmaxPabsBound}. As compared to the original structure depicted in Fig.~\ref{figmaxPabsSchematics}, the power absorbed in the \Quot{chip} was increased by a factor of~$4.52\cdot10^4$.

Comparing the objective functions in Fig.~\ref{figmaxPabsComparison}, it is clear that increasing $\Nags$ increases the number of simultaneously used local minima during the evaluation of the global step, which leads to increased diversity and, consequently, increases the chance of finding a high-quality solution. This is, of course, at the expense of computational time, where the increase is approximately linear with $\Nags$. Notice, however, that with access to the computational cluster with $C$ cores, there is almost no difference between computational times as far as $\Nags \leq C$ because all agents in the same global iteration~$j$ are evaluated simultaneously thanks to excellent scalability in parallel computing.

The optimal candidates for both runs are depicted in the bottom-right corner of Fig.~\ref{figmaxPabsComparison} which shows that the shape for~$\Nags = 192$ is more regular, consisting of three wire structures. The shape for~$\Nags = 80$ has not only poorer performance in terms of absorbed power~$P_\T{lost}$, but is also significantly more complex rendering its possible manufacturing risky.

\begin{figure}
\centering
\includegraphics[width=0.67\columnwidth, angle=90]{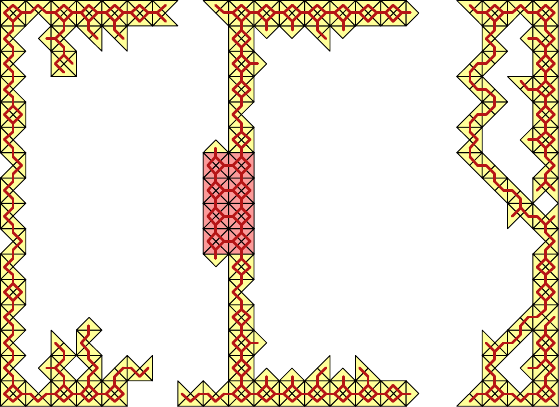}
\includegraphics[width=0.67\columnwidth, angle=90]{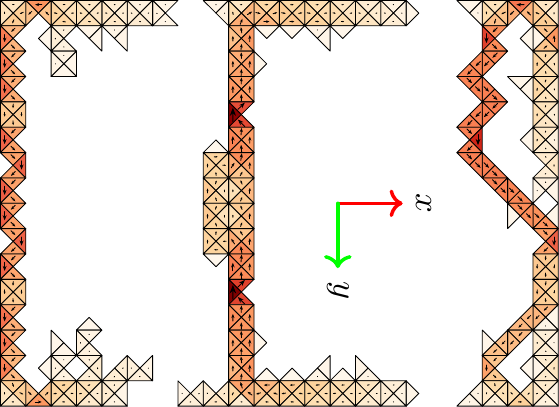}
\caption{(left) A structure optimized within a design region depicted in yellow in Fig.~\ref{figmaxPabsSchematics} for the maximization of power absorbed in a chip depicted in pink in Fig.~\ref{figmaxPabsSchematics}. The physical setting is as described in Fig.~\ref{figmaxPabsSchematics}. The number of \ac{DOF} is~$\Ndof = 905$, the number of controllable \ac{DOF} is~$\Nopt = 863$. The number of agents used was $\Nags = 192$, see the red curve in Fig.~\ref{figmaxPabsComparison}. The red connections denote enabled \ac{DOF}. (right) Current density on the optimized structure from the left pane.}
\label{figmaxPabsCurrent}
\end{figure}

\section{Discussion and Future Challenges}
\label{sec:discussion}

The method has some unique properties and offers unorthodox features which are discussed below together with a list of pros and cons.

\subsection{Properties and Features}

\subsubsection{Full-wave evaluation} The entire approach is full-wave. The numerical errors occurring during the iterative updates are negligible, typically influencing the last digit in double precision\footnote{This was verified by optimizing the shape with the proposed method first and solving the method of moments for the resulting shape directly after that. The difference in objective function was compared then for these two current densities.}.

\subsubsection{Fixed discretization grid} The procedure is based on a fixed discretization grid. It may happen that some objective function tends to thin structures, \ie{}, only one \ac{DOF} per width is chosen no matter what density of grid is used. A procedure with remeshing can be, however, applied in these cases. The objective function can also be penalized with ohmic losses, the incorporation of which always reflects physical reality better.

\subsubsection{No interpolation function} As compared to the classic topology optimization~\cite{BendsoeSigmund_TopologyOptimization, 2016_Liu_AMS}, there is no interpolation procedure which can change the actual value of an objective function when performed and which selection depends on a user.

\subsubsection{Gradient-based procedure} As opposed to the pixeling~\cite{RahmatMichielssen_ElectromagneticOptimizationByGenetirAlgorithms}, the proposed method involves the local step. Gradient-based optimization methods are typically preferred due to their rapid convergence~\cite{Strang_Applied_Math}. The gradient-based method may fail if a problem includes a significant amount of local minima with the solution space unable to be regularized~\cite{Sigmund_OnTheUselessOfNongradinetApproachesInTopoOptim}. To mitigate this, the global step is utilized. Being connected is, however, a question of whether the greedy search used in the local step is an appropriate choice as there are NP-hard problems known to be \Quot{greedy-resistant}~\cite{Bendall2006_GreedyResistantProblems}. Here we only rely on the experimental evidence presented in this paper, which suggests that the greedy search provides appropriate designs.

\subsubsection{Second topology derivative} Considering genetic algorithm as a particular choice for the global step, the mutation operator~\cite{Haupt_Werner_GeneticAlgorithmsInEM} serves as a second option. It is performed for random perturbations, but its application can be controlled by a proper choice of metaparameters.

\subsubsection{\ac{MoM} compatibility} The approach is compatible with any \ac{MoM} formulation based on piecewise basis functions and a solution found with a direct (noniterative) solver. Important differences may, however, exist depending on the type of basis functions used. For example, overlapping basis functions (\eg{}, \ac{RWG}~\cite{RaoWiltonGlisson_ElectromagneticScatteringBySurfacesOfArbitraryShape}), do not coincide one by one with discretization elements, therefore, removing one \ac{DOF} does not correspond to removing one discretization element.

\subsubsection{Compatibility with Fundamental Bounds}

The formulation within the framework of \ac{MoM} allows evaluating fundamental bound on many metrics of practical interest. These bounds then serve as stopping criteria for topology optimization, a non-convex problem in which the proximity to the global extremum is unknown.

\subsubsection{Multi-modality} The algorithm is multi-modal since, in principle, $\Nags$~local minima are found during each global step thanks to the application of the local step~\cite{Deb_MultiOOusingEA}. The number of unique local minima decreases as the optimization converges to the global minimum.

\subsubsection{Big data} A huge amount of data is gathered during the optimization, \cf{} Table~\ref{tab:QuOpt} which, together with excellent control over the entire optimization, makes this approach an ideal candidate for real-time data processing with machine learning.

\subsubsection{Flexible objective function} Topology sensitivity is evaluated as the difference between the performance of the actual current~$\Iv \left(\gene_i\right)$ and the performance of the current flowing on the perturbed structure~$\Iv \left(\gene_{i+1}\right)$. As such, the objective function can be easily extended towards:
\begin{itemize}
    \item an objective function evaluated only within a sub-region (see Section~\ref{sec:maxPabs}),
    \item an objective function taking into account multiple frequency samples,
    \item an objective function involving port quantities defined for multiple ports.
\end{itemize}

\subsubsection{Multi-objective optimization} Multi-objective optimization can easily be performed with the proposed memetics by scalarization, see Section.~\ref{sec:TradeOff}. Only convex trade-offs can, however, be found with a weighted sum, such as in~\eqref{eq:Qfact} or~\eqref{eq:QminZin} where more advanced techniques would have to be adopted for complex Pareto frontiers, \eg{}, a rotated weighted metric method~\cite{Deb_MultiOOusingEA}.

\subsection{Generalization of Variable Space}

\subsubsection{Feeding synthesis} The optimal placement of the feeder(s), including optimal amplitudes and phases, can be incorporated into shape optimization. For one feeder, it is, technically, not needed\footnote{It is assumed that the optimal shape is adjusted according to the initial placement of the feeder.}. For more feeders, it is an extra combinatorial problem which can run over or run together with the shape optimization.

\subsubsection{Optimization variables} It is possible to generalize the entire algorithm so that discretization elements are dealt with and not \ac{DOF} (basis functions). One possibility is to define the smallest blocks manually, and optimize over them only. This slows down the evaluation since vectorization is not as efficient because block inversion has to be applied. On the other hand, the number of unknowns is reduced depending on the number and size of the optimized blocks. The initial tests indicate that computational time is comparable or lower than with \ac{DOF} as the optimization unknowns.

\subsubsection{Optimization domain} Topology sensitivity can be evaluated only for removals or additions offering a substantial speed-up. Alternatively, as often required in practice, only a part of the structure can be optimized, see Fig.~\ref{figmaxPabsSchematics}, reducing computational time proportionally to the size of the optimized region. Another possibility is to investigate only \ac{DOF}, which separates two regions filled by different materials (the number of unknowns is reduced approximately by one order).

\subsubsection{Multi-material representation} The two-state, binary, material representation (PEC $\times$ vacuum) can easily be generalized to a multi-state optimization at the cost of the linear increase of variable space size.

\subsection{Deficiencies and Possible Remedies}

\subsubsection{Slow evolution} 
\begin{figure}
\centering
\includegraphics[width=\columnwidth]{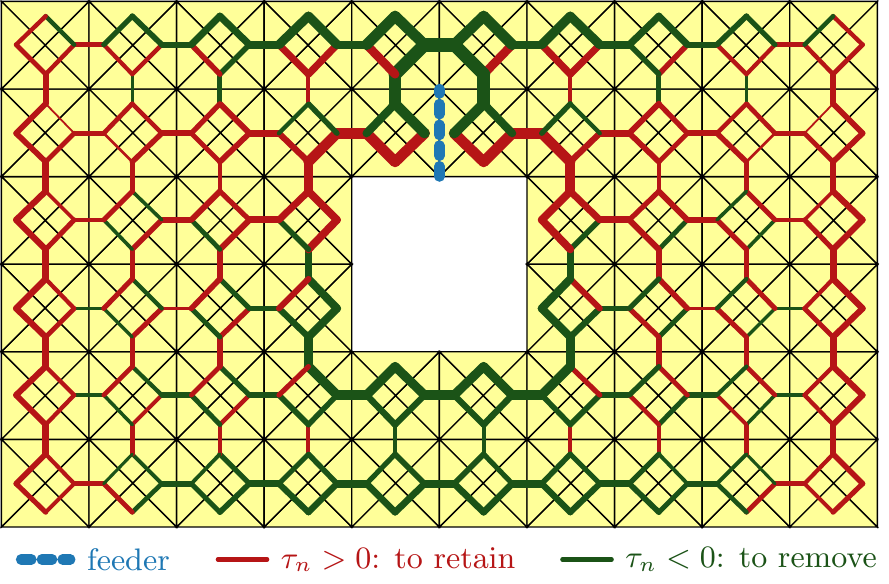}
\caption{Topology sensitivity map of a rectangular PEC plate of electrical size~$ka = 0.5$ and a hole in the middle. The objective function evaluates Q-factor~\eqref{eq:Qfact}. The delta gap feeder is depicted by the thick blue line. The sensitivites~$\tau_n$ are evaluated for all $\Ndof$ \ac{DOF}. The negative sensitivities, $\tau_n < 0$, are highlighted by the green lines, the positive sensitivies, $\tau_n > 0$, by the red lines. The thicker the line, the higher the absolute value of the sensitivity.}
\label{figTSeval}
\end{figure}

Shape modification is performed via the removal or addition of a \ac{DOF} decreasing the value of the objective function the most. Following this greedy search, the \ac{DOF} corresponding to the thickest green line in Fig.~\ref{figTSeval}, \ie{}, the connections located to the left and on the top right section of the feeder are to be removed in order to eliminate the short-circuit. Such an approach is, however, relatively slow\footnote{As compared to the adjoint formulation of topology optimization~\cite{BendsoeSigmund_TopologyOptimization} where all \ac{DOF} are updated at once.}, requiring in many cases~$\OP{O}(\Ndof)$ local updates. The possible remedy might be to average the sensitivity in a small area and update all \ac{DOF} inside at once. 

\subsubsection{Algorithm complexity} Numerical complexity is unpleasant as it  favors small to mid-size structures. Realistically, with the current implementation and state-of-the-art hardware, thousands of unknowns are possibly optimized in hours. Billions of \ac{DOF} are, however, common in FEM~\cite{Aage_etal_TopologyOptim_Nature2017}, but this is because of different algebraic properties of the stiffness matrix as compared to the impedance matrix.

\subsubsection{Shape irregularity} 
\begin{figure}
\centering
\includegraphics[width=\columnwidth]{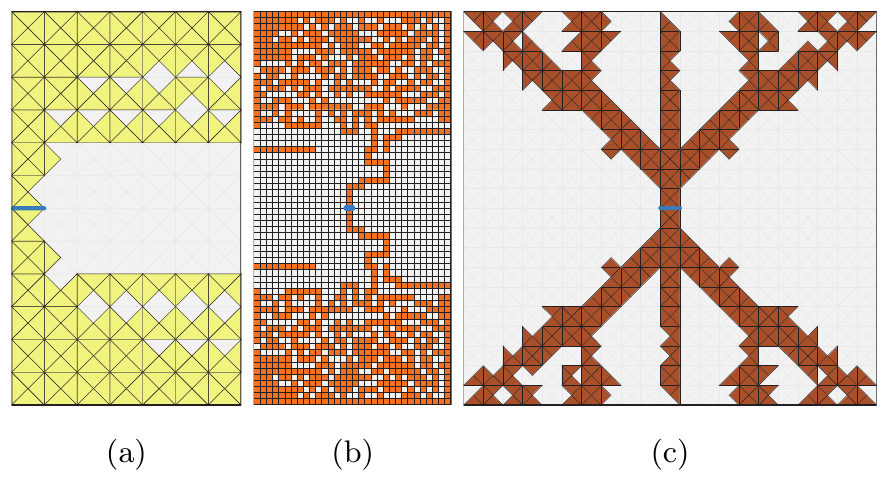}
\caption{Comparison of final structures optimized with various shape optimization techniques. (a) The method proposed in this work, Q-factor minimization. (b) Pixeling based on the genetic algorithm, Q-factor minimization~\cite{CismasuGustafsson_FBWbySimpleFreuqSimulation}. (c) Adjoint formulation of topology optimization, total efficiency minimization~\cite{2016_Liu_AMS}. All techniques utilize method of moments to evaluate the objective function. The orange dashed line delimits the design region.}
\label{figRegularity}
\end{figure}

The method might produce highly irregular optimized shapes. This is a common issue of shape optimization, no matter what technique is used, see Fig.~\ref{figRegularity}, comparing representatives of this method, pixeling based on genetic algorithm~\cite{RahmatMichielssen_ElectromagneticOptimizationByGenetirAlgorithms}, and gradient-based topology optimization~\cite{2016_Liu_AMS}. All shapes are relatively irregular, \ie{}, this property is common to shape optimization routines. Hence, the next step could be to find a way how to regularize the shapes.

\section{Conclusion}
\label{sec:concl}

A novel memetic procedure for optimizing electromagnetic devices was presented, consisting of local and global steps to mitigate the disadvantages of both approaches when used alone. The fixed discretization grid is assumed, allowing the method of moments system matrix to be inverted only once while storing it and other required matrices in the computer's memory. The iterative full-wave evaluation of all the smallest topology perturbations is done using an inversion-free (Sherman-Morrison-Woodbury) formula and an exact-reanalysis-based procedure. The local step provides gradient-type information about the topology, and local updates are performed in a greedy sense, \ie{} with the part of the optimized shape enhancing the performance the most being updated in each local iteration. The global step used in this paper is based on a genetic algorithm which maintains diversity and increases the chance of the algorithm's convergence towards a global minimum. Genetic algorithms match well with the discrete form of the optimization problem and, as a result, the memetic procedure is robust, fast, and versatile.

The procedure has many advantages and unique features, such as being able to find minima close to the fundamental bounds, implying that the local minima are close to the global one. The number of full-wave evaluated shapes spans from millions for small problems to billions for larger problems with thousands of unknowns. The optimization approach is multi-modal and capable of finding many local minima at once during one global step. As a result, the global scheme optimizes shapes only in a significantly reduced solution space containing only local minima. Since the optimized shape is known at every step, a multi-port, multi-frequency, or multi-material optimization can be performed only at the expense of the linear increase of computational time.

Examples shown in this part proved the efficiency of the method. For the first time, the performance of optimized structures has been directly compared with the fundamental bounds on the optimized metrics. This practice provides an ultimate measure of optimization efficiency, but also naturally scales optimized metrics and establishes the straightforward terminal criterion for optimization. In many cases, the performance of the resulting shapes closely follows the bounds. The excellent performance of the algorithm demonstrates its wide applicability in various inverse design problems in electromagnetics.

There are still many future challenges that could broaden the usage of the proposed method to make it even more efficient, such as including lumped elements and their synthesis. It can be shown that this technique is similar to adjoint formulation over gray-scaled material. Another possibility is to introduce a set of geometrical operators, representing geometrical metrics, such as the area spanned by the material, or the curvature of the shape. Having complete freedom in formulating the objective function, these operators may improve the regularity of shapes and eliminate manufacturing difficulties. Technical, and still extremely relevant, improvements include GPU implementation or a detailed study of metaparameters used for optimization settings and their tuning at the beginning or throughout an optimization. Another way to accelerate the evaluation is to use an adaptive scheme with a successively refined discretization grid to impose the optimal structure from a coarse grid as the initial shape for the finer grid.

\appendices

\section{Fundamental Bound on Antenna Q-factor}
\label{app:Appendix0}

The lower bound on radiation Q-factor~$Q_\T{lb}$ is found via a \ac{QCQP}
\begin{equation}
\begin{split}
\underset{\Iv}{\T{minimize}} \quad & \Iv^\herm \Wm \Iv \\
\T{subject\, to} \quad & \Iv^\herm \Rmvac \Iv = 1 \\
& \Iv^\herm \Xm \Iv = 0,
\end{split}
\label{eq:Qmin1}
\end{equation}
where~$\Rmvac$ is radiation and~$\Xm$ is the reactance part of impedance matrix~$\Zm$, and $\Wm$~is the stored energy matrix. The problem is recast into its dual form and solved via a generalized eigenvalue problem as described in~\cite{Liska_etal_FundamentalBoundsEvaluation}. 

Q-factor~$Q_\T{lb}^\T{TM}$ is a tighter bound for all cases when only the TM modes are involved. This includes, \eg{}, planar antennas with a discrete feeder~\cite{Best_ElectricallySmallResonantPlanarAntennas, Capek_etal_2019_OptimalPlanarElectricDipoleAntennas}, \ie{}, the case studied in Section~\ref{sec:QminU}. The formula~\eqref{eq:Qmin1} still applies with a change~$\M{R}_0 \to \M{R}_0^\T{TM}$, where matrix~$\M{R}_0^\T{TM}$ is evaluated as
$\M{R}_0^\T{TM} = \left( \M{U}_1^\T{TM} \right)^\herm \M{U}_1^\T{TM}$
with~$\M{U}_1^\T{TM}$ being a projection matrix from the basis of transverse-magnetic (TM) spherical vector waves into an \ac{MoM} basis~\cite{2022_Losenicky_TAP}.

\section{Fundamental Bound on Realized Gain With Prescribed Input Impedance}
\label{app:AppendixA}

The upper bound on realized gain was evaluated using~\eqref{eq:GmaxZin} with current vector~$\M{I}$ being the solution to
\begin{equation}
\begin{split}
\underset{\Iv}{\T{minimize}} \quad & - \Iv^\herm \M{F}^\herm (\UV{d}, \UV{e}) \M{F} (\UV{d}, \UV{e})  \Iv \\
\T{subject\, to} \quad & \Iv^\herm \M{Z} \Iv = \Iv^\herm \M{V} \\
& \M{V}^\herm \Iv = \dfrac{\left| V_\T{in} \right|^2}{\Zin^0},
\end{split}
\label{eq:Gain}
\end{equation}
where the last affine constraint enforcing matching~$\varGamma = 0$, see~\eqref{eq:GmaxZin}, is removed from the optimization using the affine transformation described in~\cite{Liska_etal_FundamentalBoundsEvaluation}. Voltage~$V_\T{in}$, imposed on the delta-gap source together with a real impedance~$\Zin^0$, is related to input power and matrix~$\M{F} (\UV{d}, \UV{e})$ gives the electric far field.

\section{Fundamental Bound on Absorbed Power in an Uncontrollable Subregion}
\label{app:AppendixB}

This fundamental bound assumes partitioning
\begin{equation}
\label{eq:cu1}
\mqty [
\M{Z}_\T{cc} & \M{Z}_\T{cu}\\
\M{Z}_\T{uc} & \M{Z}_\T{uu}
]
\mqty[\M{I}_\T{c}\\ \M{I}_\T{u}] 
= 
\mqty[ \M{V}_\T{c}\\ \M{V}_\T{u} ]
\end{equation}
into ``controllable'' (index c) and ``uncontrollable'' (index u) subregions. The QCQP leading to the desired optimal current reads
\begin{equation}
\begin{aligned}
	\min \limits_{\M{I}_\T{c}} \quad & - \M{I}^\herm \mqty [
\M{0}_\T{cc} & \M{0}_\T{cu}\\
\M{0}_\T{uc} & \M{R}_{\rho,\T{uu}}
] \M{I} \\
	\T{s.t.} \quad &\M{I}^\herm \M{Z} \M{I} = \M{I}^\herm \M{V} \\
	&\mqty [\M{Z}_\T{uc} & \M{Z}_\T{uu}]\M{I}= \M{V}_\T{u},
\end{aligned}  
\label{eq:Pauu_opt}
\end{equation}
where the first constraint enforces the conservation of complex power~\cite{2020_Gustafsson_NJP} and where the last affine constraint enforcing the bottom row of the partitioned system~\eqref{eq:cu1} is, as in Appendix~\ref{app:AppendixA}, removed from the optimization using the affine transformation described in~\cite{Liska_etal_FundamentalBoundsEvaluation}.
 \section*{Acknowledgement}
The access to the computational infrastructure of the OP VVV funded project CZ.02.1.01/0.0/0.0/16\_019/0000765 ``Research Center for Informatics'' is gratefully acknowledged.


\begin{IEEEbiography}[{\includegraphics[width=1in,height=1.25in,clip,keepaspectratio]{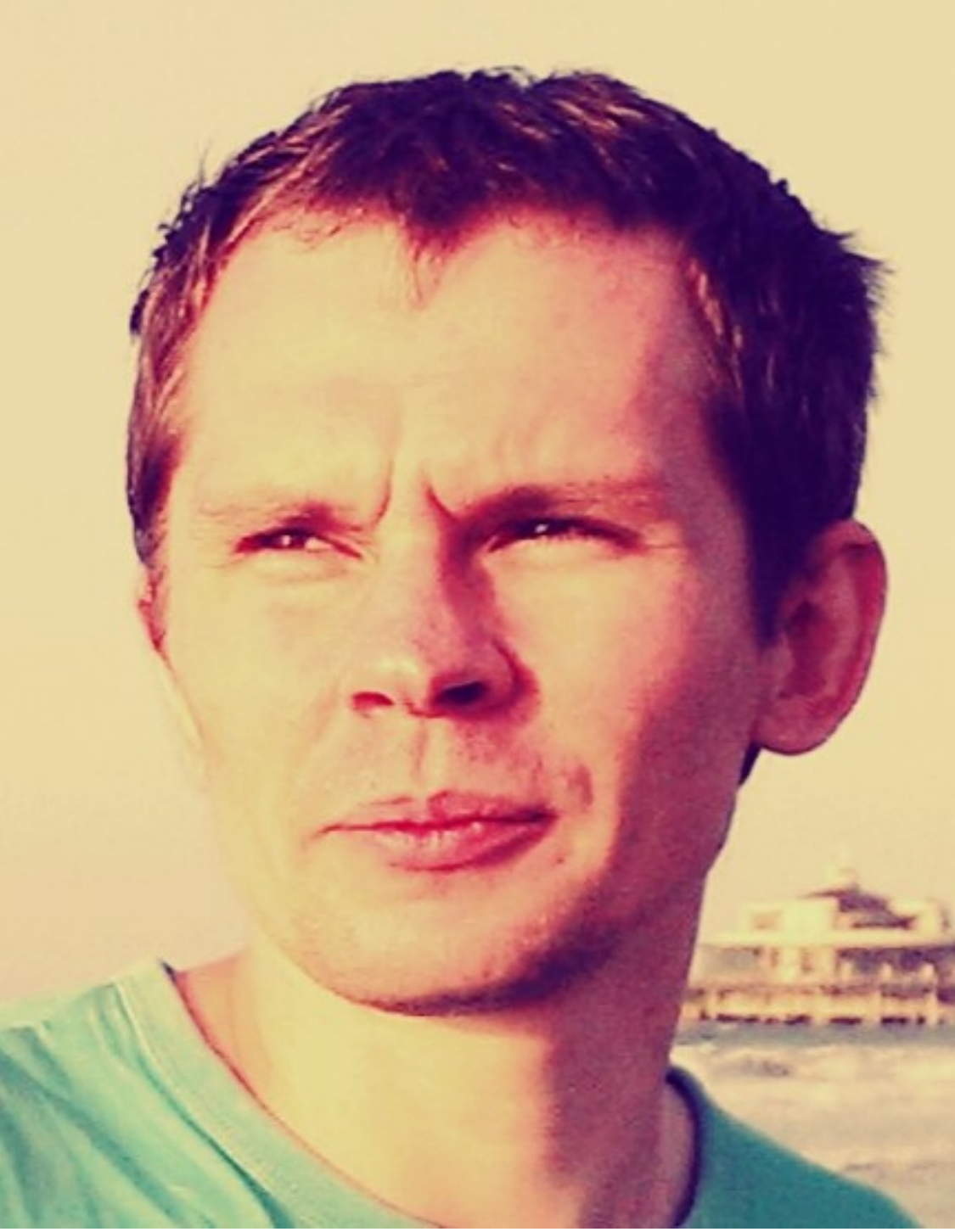}}]{Miloslav Capek}
(M'14, SM'17) received the M.Sc. degree in Electrical Engineering 2009, the Ph.D. degree in 2014, and was appointed Associate Professor in 2017, all from the Czech Technical University in Prague, Czech Republic.
	
He leads the development of the AToM (Antenna Toolbox for Matlab) package. His research interests are in the area of electromagnetic theory, electrically small antennas, numerical techniques, fractal geometry, and optimization. He authored or co-authored over 120~journal and conference papers.

Dr. Capek is Associate Editor of IET Microwaves, Antennas \& Propagation. He was a regional delegate of EurAAP between 2015 and 2020. He received the IEEE Antennas and Propagation Edward E. Altshuler Prize Paper Award 2023.
\end{IEEEbiography}

\begin{IEEEbiography}[{\includegraphics[width=1in,height=1.25in,clip,keepaspectratio]{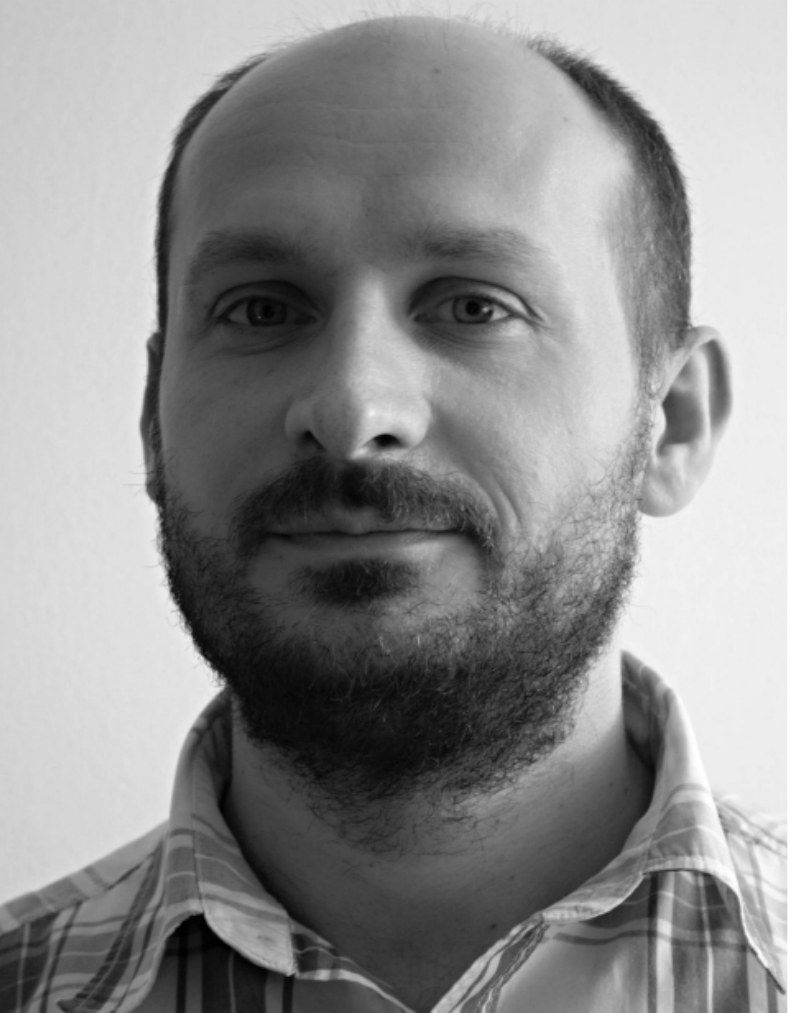}}]{Lukas Jelinek}
received his Ph.D. degree from the Czech Technical University in Prague, Czech Republic, in 2006. In 2015 he was appointed Associate Professor at the Department of Electromagnetic Field at the same university.

His research interests include wave propagation in complex media, electromagnetic field theory, metamaterials, numerical techniques, and optimization.
\end{IEEEbiography}

\begin{IEEEbiography}[{\includegraphics[width=1in,height=1.25in,clip,keepaspectratio]{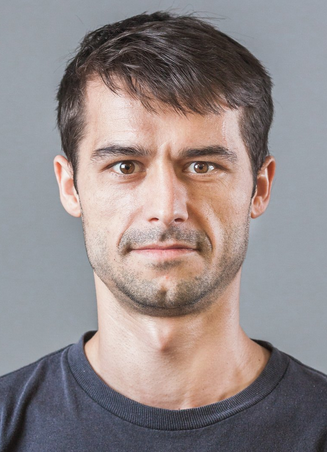}}]{Petr Kadlec}
(M'13) received the Ph.D. degree in electrical engineering from the Brno University of Technology (BUT), Brno, Czech Republic, in 2012. He is currently an Associate Professor with the Department of Radioelectronics, BUT. His research interests include global optimization methods and computational methods in electromagnetics. He is a leading developer of the FOPS (Fast Optimization ProcedureS) MATLAB software package. 
\end{IEEEbiography}

\begin{IEEEbiography}[{\includegraphics[width=1in,height=1.25in,clip,keepaspectratio]{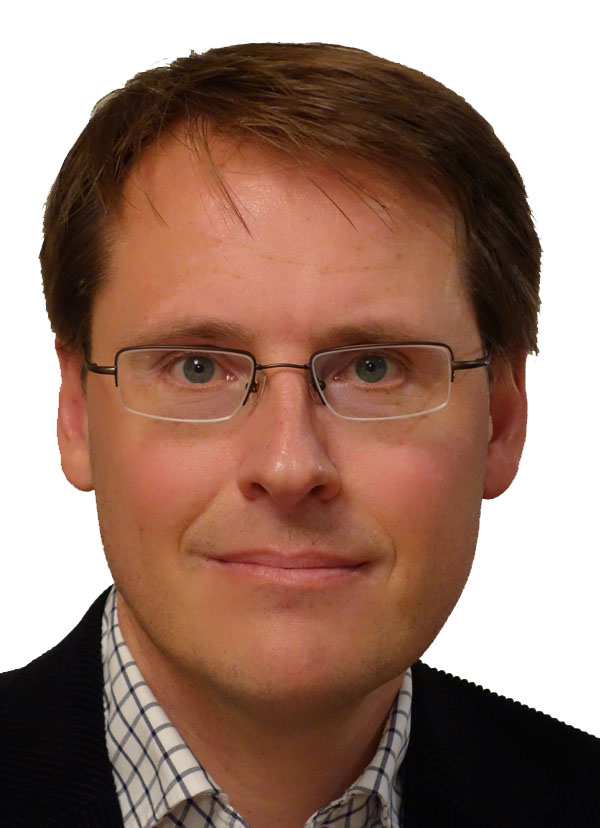}}]{Mats Gustafsson}
received the M.Sc. degree in Engineering Physics 1994, the Ph.D. degree in Electromagnetic Theory 2000, was appointed Docent 2005, and Professor of Electromagnetic Theory 2011, all from Lund University, Sweden.

He co-founded the company Phase holographic imaging AB in 2004. His research interests are in scattering and antenna theory and inverse scattering and imaging. He has written over 100 peer reviewed journal papers and over 100 conference papers. Prof. Gustafsson received the IEEE Schelkunoff Transactions Prize Paper Award 2010, the IEEE Uslenghi Letters Prize Paper Award 2019, and best paper awards at EuCAP 2007 and 2013. He served as an IEEE AP-S Distinguished Lecturer for 2013-15.
\end{IEEEbiography}

\end{document}